%% file: berard-saearp-mat-contemp-2008.tex
\documentclass[10pt,twoside]{article}
\usepackage[frenchb,english]{babel}
\usepackage{version}

\usepackage{graphicx}

\includeversion{pb1-figs} 
\excludeversion{file-name}
\begin{file-name}
\usepackage{fink}
\end{file-name}\bigskip

\input{macros}

\begin{document}

\title{\vspace{-3cm}\hspace{5.0cm} \phantom{xxxxxx} {\scriptsize to appear in \em Matem\'atica
Contempor\^anea 34 (2008)}
\vspace*{1cm} \\
\centering  \large \bf EXAMPLES OF  $H$-HYPERSURFACES IN $\HH^n
\times \R$  \\AND GEOMETRIC APPLICATIONS \footnotetext{ {Mathematics
Subject Classification (2000):~} 53C42, 53C21
\newline
{\em Key words:~} Constant mean curvature, Rotation and
Translation hypersurfaces, Vertical graphs, Dirichlet problem.
}}
\author{\large \bf P. B\'{e}rard \and \large \bf  R. Sa Earp}
\date{}
\maketitle

\begin{center}
{\small \it
 Dedicated to Professor Manfredo do Carmo on the
occasion of his 80$^{\text{th}}$ birthday}
\end{center}

\thispagestyle{empty}
\pagestyle{myheadings}
\markboth{P. BERARD \quad R. SA EARP}
  {EXAMPLES OF  $H$-HYPERSURFACES IN $\HH^n \times \R$  AND GEOMETRIC APPLICATIONS }
\setlength{\baselineskip}{6mm}

\pagestyle{myheadings}
\thispagestyle{empty}

\begin{abstract}
\noi In this paper we describe all rotation $H$-hypersurfaces in
$\HH^n \times \R$ and use them as barriers to prove existence and
characterization of certain vertical $H$-graphs and to give symmetry
and uniqueness results for compact $H$-hypersurfaces whose boundary
is one or two parallel submanifolds in slices. We also describe
examples of translation $H$-hypersurfaces in $\HH^n \times \R$. For
$n>2$ we obtain a complete embedded translation hypersurface
generated by a compact, simple, strictly convex curve. When $0 < H <
\frac{n-1}{n}$ we obtain a complete non-entire vertical graph over
the non-mean convex domain bounded by an equidistant hypersurface
taking infinite boundary value data and infinite asymptotic boundary
value data.
\end{abstract}

\section{Introduction}\label{S-intro}


Rotation and translation surfaces with constant mean curvature in
$\HH^2 \times \R$ have been studied in details in \cite{ST05, Sa08,
ST08} together with applications. We have studied rotation and
translation minimal hypersurfaces with applications in
\cite{BS08a}.\bigskip

In this paper, we consider constant non-zero mean curvature
hypersurfaces in $\HH^n \times \R$.\bigskip

We consider rotation $H$-hypersurfaces in Section
\ref{SS-dim3-Hrot}. For $H > \frac{n-1}{n}$, we find the constant
mean curvature sphere-like hypersurfaces obtained in \cite{HH89} and
the Delaunay-like hypersurfaces obtained in \cite{PR99}. When $0< H
\le \frac{n-1}{n}$, we obtain complete simply-connected
hypersurfaces $\cS_H$ which are entire vertical graphs above
$\HH^n$, as well as some complete embedded or complete immersed
cylinders which are bi-graphs (Theorems~\ref{T-h3-r0} and
\ref{T-h3-rp}). When $H=\frac{n-1}{n}$, the asymptotic behaviour of
the height function of these hypersurfaces is exponential, and it
only depends on the dimension when $n \ge 3$. In
Section~\ref{S-appli}, we give geometric applications using the
simply-connected rotation $H$-hypersurfaces $\cS_H$ ($0 < H \le
\frac{n-1}{n}$) mentioned above as barriers. We give existence and
characterization of vertical $H$-graphs ($0 < H \le \frac{n-1}{n}$)
over appropriate bounded domains (Proposition \ref{P-appl-2}) as
well as symmetry and uniqueness results for compact hypersurfaces
whose boundary is one or two parallel submanifolds in slices
(Theorems \ref{T-appl-6} and \ref{T-appl-7}). These results
generalize the $2$-dimensional results obtained previously in
\cite{NSST08}.\bigskip

We treat translation $H$-hypersurfaces in Section
\ref{SS-dim3-Htransl} (Theorem \ref{T-tra-1}). When $n\ge 3$ and
$H=\frac{n-1}{n}$, we in particular find a complete embedded
hypersurface generated by a compact, simple, strictly convex
curve.\bigskip

When $0 < H <\frac{n-1}{n}$, we obtain a complete non-entire
vertical graph over the non-mean convex domain bounded by an
equidistant hypersurface $\Gamma$. This graph takes infinite
boundary value data on $\Gamma$ and it has infinite asymptotic
boundary value data. \bigskip

The authors would like to thank the Mathematics Department of
PUC-Rio (PB) and the Institut Fourier -- Universit\'{e} Joseph Fourier
(RSA) for their hospitality. They gratefully acknowledge the
financial support of CNPq, FAPERJ (in particular \emph{Pronex} and
\emph{Cientistas do nosso Estado}), Acordo Brasil - Fran\c{c}a,
Universit\'{e} Joseph Fourier and R\'{e}gion Rh\^{o}ne-Alpes.\bigskip

\bigskip
\section{Examples of $H$-hypersurfaces in $\HH^n \times \R$}\label{S-examples}


We consider the ball model for the hyperbolic space $\HH^n$,

$$\B := \ens{(x_1, \ldots , x_n) \in \R^n}{x_1^2 + \cdots + x_n^2 < 1},$$

with the hyperbolic metric $g_{\B}$,

$$g_{\B} := 4 \big( 1 - (x_1^2 + \cdots + x_n^2) \big)^{-2} \big( dx_1^2 +
\cdots + dx_n^2\big),$$

and the product metric

$$\gh = g_{\B} + dt^2$$

on $\HH^n \times \R$. \bigskip


\subsection{Rotation $H$-hypersurfaces in $\HH^n \times \R$}
\label{SS-dim3-Hrot} \bigskip

The mean curvature equation for rotation hypersurfaces,

\begin{equation*}\label{E-h3-rot1a}
n H(\rho) \sinh^{n-1}(\rho) = \partial_{\rho} \Big(
\sinh^{n-1}(\rho) \dot{\lambda}(\rho) (1 +
\dot{\lambda}^2(\rho))^{-1/2}\Big)
\end{equation*}\medskip

can be established using the flux formula, see Appendix
\ref{S-vflux}. We consider rotation hypersurfaces about $\{0\}
\times \R$, where $\rho$ denotes the hyperbolic distance to the axis
and the mean curvature is taken with respect to the unit normal
pointing upwards. \bigskip

Minimal rotation hypersurfaces in $\HH^n \times \R$ have been
studied in \cite{ST05} in dimension $2$ and in \cite{BS08a} in
higher dimensions. In this Section we consider the case in which $H$
is a non-zero constant. We may assume that $H$ is positive.\bigskip

Integrating the above differential equation, we obtain the equation
for the generating curves of rotation $H$-hypersurfaces in $\HH^n
\times \R$,

\begin{equation}\label{E-h3-rot1}
\dot{\lambda}(\rho) \big( 1+ \dot{\lambda}^2(\rho)\big)^{-1/2}
\sinh^{n-1}(\rho) = nH \int_0^{\rho} \sinh^{n-1}(t) \, dt + d
\end{equation}

for $H > 0$ and for some constant $d$. \bigskip

This equation has been studied in \cite{NSST08, ST05} in dimension
$2$ (with a different constant $d$). \bigskip

\noi \textbf{Notations.}~ For later purposes we introduce some
notations. \bigskip

\noib For $m \ge 0$, we define the function $I_m(t)$ by
\begin{equation}\label{E-h3-rot2}
I_m(t) := \int_0^{t} \sinh^m(r) \, dr.
\end{equation} \bigskip

\noib For $H>0$ and $d \in \R$, we define the functions,

\begin{equation}\label{E-h3-4}
\left\{%
\begin{array}{lll}
    M_{H,d}(t) &:=& \sinh^{n-1}(t) - nH I_{n-1}(t) - d,\\
    P_{H,d}(t) &:=& \sinh^{n-1}(t) + nH I_{n-1}(t) + d,\\
    Q_{H,d}(t) &:=& \big[nH I_{n-1}(t) + d\big] \big[M_{H,d}(t) \,
    P_{H,d}(t)\big]^{-1/2}, \\
    && \text{when the square root exists.}\\
\end{array}%
\right.
\end{equation}\bigskip

We see from (\ref{E-h3-rot1}) that $\dot{\lambda}(t)$ has the sign
of $n H I_{n-1}(t) +d$. It follows that $\lambda$ is given, up to an
additive constant, by

$$\lambda_{H,d}(\rho) = \int_{\rho_0}^{\rho} \frac{nH I_{n-1}(t) +
d}{\sqrt{\sinh^{2n-2}(t) - \big(nH I_{n-1}(t) + d\big)^2}}\, dt$$

or, with the above notations,

\begin{equation}\label{E-h3-rot3}
\lambda_{H,d}(\rho) = \int_{\rho_0}^{\rho} \frac{nH
I_{n-1}(t)+d}{\sqrt{M_{H,d}(t) P_{H,d}(t)}}\, dt =
\int_{\rho_0}^{\rho} Q_{H,d}(t) \, dt
\end{equation}

where the integration interval $[\rho_0,\rho]$ is contained in the
interval in which the square-root exists. The existence and
behaviour of the function $\lambda_{H,d}$ depend on the signs of the
functions $nH I_{n-1}(t) + d$, $M_{H,d}(t)$ and $P_{H,d}(t)$.
\bigskip

Up to vertical translations, the rotation hypersurfaces about the
axis $\{0\} \times \R$, with constant mean curvature $H>0$ with
respect to the unit normal pointing upwards, can be classified
according to the sign of $H - \frac{n-1}{n}$ and to the sign of $d$.
We state three theorems depending on the value of $H$. \bigskip

\newpage

\begin{thm}[Rotation $H$-hypersurfaces with $H=\frac{n-1}{n}$]\label{T-h3-r0}
$ $
    \begin{enumerate}
    \item When $d=0$, the hypersurface $\cS_{\frac{n-1}{n}}$ is a
    simply-connected entire vertical graph above $\HH^n \times \{0\}$,
    tangent to the slice at $0$, generated by a strictly convex
    curve. The height function $\lambda (\rho)$ on $\cS_{\frac{n-1}{n}}$
    grows exponentially.

    \item When $d>0$, the hypersurface $\cC_{\frac{n-1}{n}}$ is a complete
    embedded cylinder, symmetric with respect to the slice $\HH^n \times \{0\}$.
    The parts $\cC_{\frac{n-1}{n}}^{\pm} := \cC_{\frac{n-1}{n}} \cap \HH^n
    \times \R_{\pm}$ are vertical graphs above the exterior of a ball
    $B(0,a)$, for some constant $a > 0$ depending on $d$.
    The height function $\lambda (\rho)$ on $\cC_{\frac{n-1}{n}}^{\pm}$ grows
    exponentially. When $n=2$, the solution exists when $0< d<1$ only.

    \item  When $d<0$, the hypersurface $\cD_{\frac{n-1}{n}}$ is complete and
    symmetric with respect to the slice $\HH^n \times \{0\}$. It has
    self-intersections along a sphere in $\HH^n \times \{0\}$. The parts
    $\cD_{\frac{n-1}{n}}^{\pm} := \cD_{\frac{n-1}{n}} \cap \HH^n
    \times \R_{\pm}$ are vertical graphs above the exterior of a ball
    $B(0,a)$, for some constant $a > 0$ depending on $d$. The height
    function $\lambda (\rho)$ on $\cD_{\frac{n-1}{n}}^{\pm}$ grows exponentially.
    \end{enumerate}

    The asymptotic behaviour of the height function when $\rho$ tends to
    infinity is as follows.
$$
\left\{%
\begin{array}{l}
\text{For~ } n=2, ~\lambda (\rho) \sim
    \frac{e^{\rho/2}}{\sqrt{1 - d}}.\\[8pt]
\text{For~ } n=3, ~\lambda (\rho) \sim
    \frac{1}{2 \sqrt{2}} \int^{\rho} \frac{e^{t}}{\sqrt{t}}\, dt.\\[8pt]
\text{For~ } n \ge 4, ~\lambda (\rho) \sim a(n) e^{b(n)t}, \text{
~for
some positive constants~ } a(n), b(n).\\
\end{array}%
\right.
$$
\end{thm}

The generating curves are obtained by symmetries from the curves
$(=)$ (standing for $H = \frac{n-1}{n}$) which appear in
Figures~\ref{F-rot-1}-\ref{F-rot-3}. \bigskip

\begin{pb1-figs}
\begin{figure}[h!]
\begin{center}
\begin{minipage}[c]{6.5cm}
    \includegraphics[width=6.5cm]{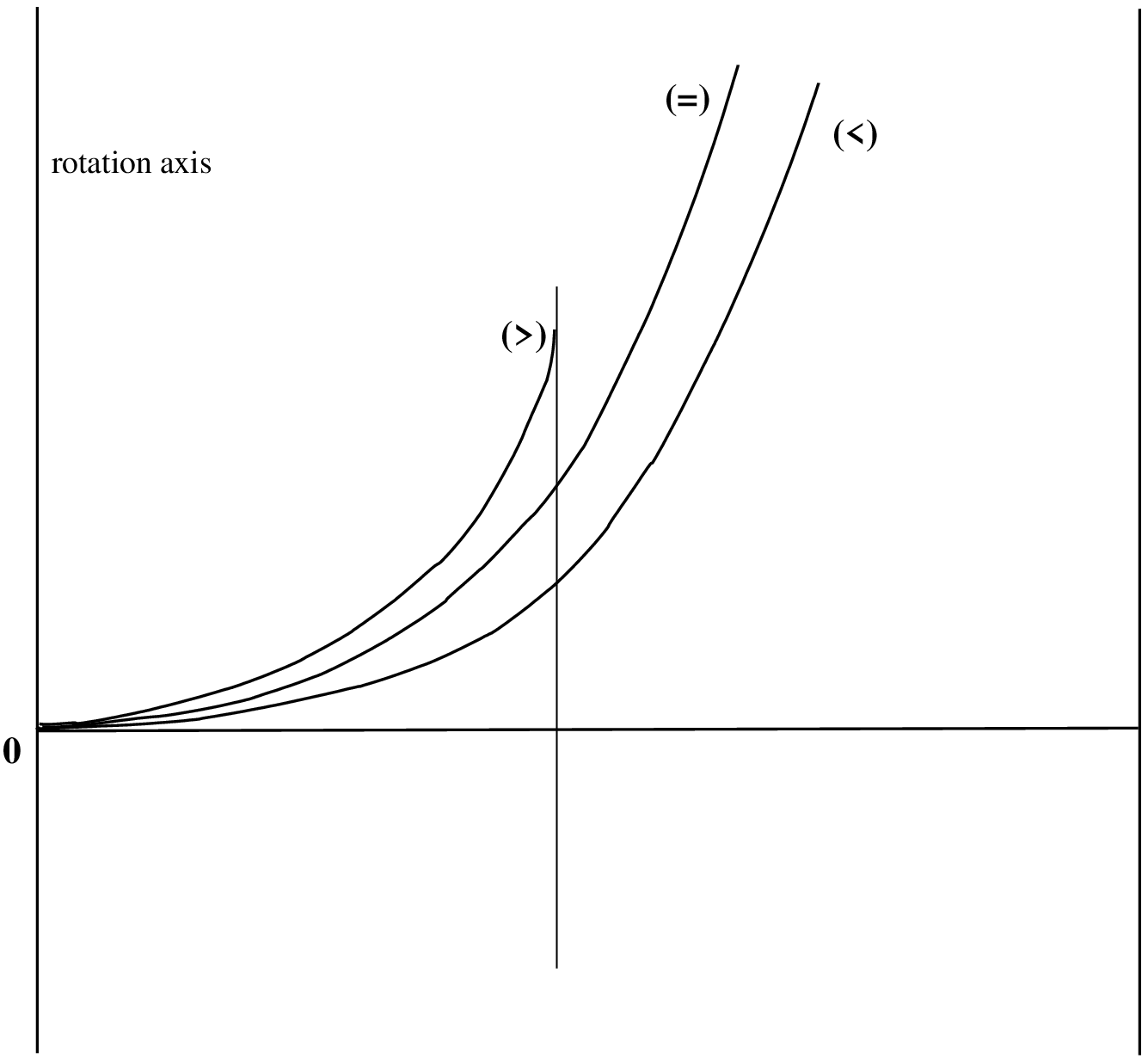}
    \caption[Case $d=0$]{Case $d=0$}
    \label{F-rot-1}
\end{minipage}\hfill
\begin{minipage}[c]{6.5cm}
    \includegraphics[width=6.5cm]{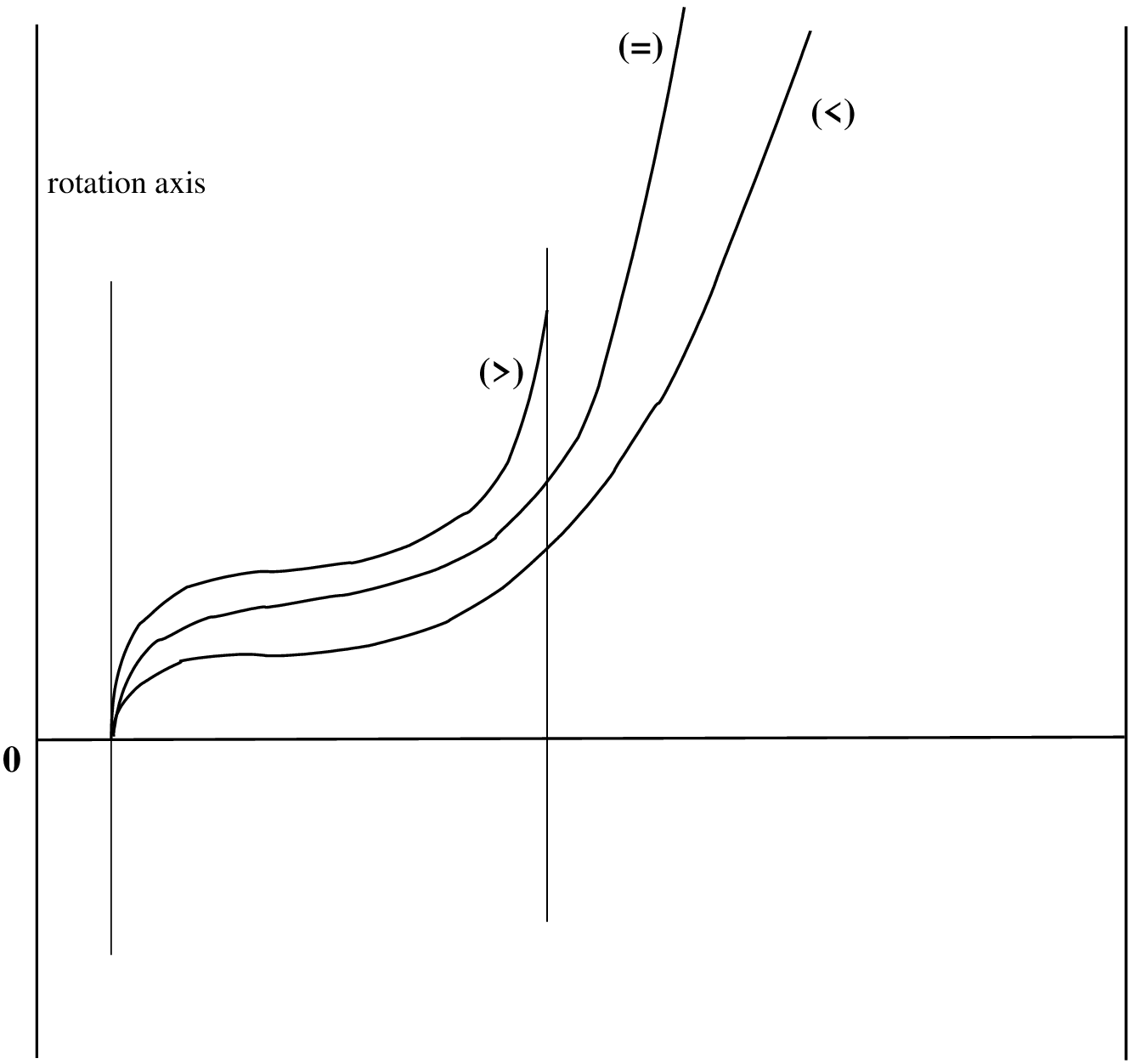}
    \caption[Case $d>0$]{Case $d>0$}
    \label{F-rot-2}
\end{minipage}\hfill
\end{center}
\end{figure}\bigskip
\end{pb1-figs}

\begin{pb1-figs}
\begin{figure}[h!]
\begin{center}
\begin{minipage}[c]{6.5cm}
    \includegraphics[width=6.5cm]{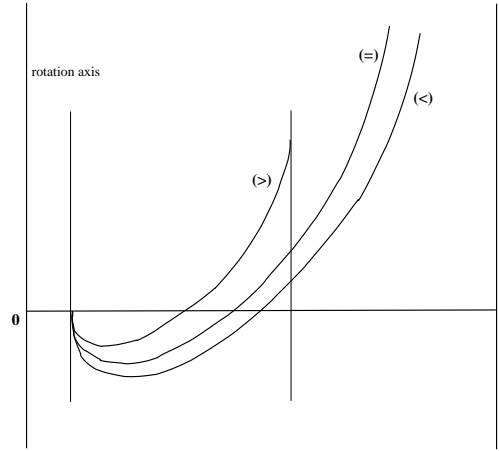}
    \caption[Case $d<0$]{Case $d<0$}
    \label{F-rot-3}
\end{minipage}\hfill
\end{center}
\end{figure}\bigskip
\end{pb1-figs}

\textbf{Remark.}~ When $n=2$ the asymptotic growth depends on the
value of the integration contant $d$. \bigskip

\begin{thm}[Rotation $H$-hypersurfaces with $0 < H < \frac{n-1}{n}$]\label{T-h3-rp}
$ $
    \begin{enumerate}
    \item When $d=0$, the hypersurface $\cS_{H}$ is a
    simply-connected entire vertical graph above $\HH^n \times \{0\}$,
    tangent to the slice at $0$, generated by a strictly convex
    curve. The height function $\lambda (\rho)$ on $\cS_{H}$ grows linearly.

    \item When $d>0$, the hypersurface $\cC_{H}$ is a complete embedded cylinder,
    symmetric with respect to the slice $\HH^n \times \{0\}$. The parts
    $\cC_{H}^{\pm} := \cC_{H} \cap \HH^n \times \R_{\pm}$ are vertical graphs above
    the exterior of a ball $B(0,a)$, for some constant $a > 0$ depending on
    $H$ and $d$. The height function $\lambda (\rho)$ on $\cC_{H}^{\pm}$ grows linearly.

    \item When $d<0$, the hypersurface $\cD_{H}$ is complete and symmetric with
    respect to the slice $\HH^n \times \{0\}$. It has self-intersections along a
    sphere in $\HH^n \times \{0\}$. The parts $\cD_{H}^{\pm} := \cD_{H} \cap
    \HH^n \times \R_{\pm}$ are vertical graphs above
    the exterior of a ball $B(0,a)$, for some constant $a > 0$ depending on $H$
    and $d$. The height function $\lambda(\rho)$ on $\cD_{H}^{\pm}$ grows linearly.
    \end{enumerate}

The asymptotic behaviour of the height function when $\rho$ tends to
infinity is given by
$$\lambda (\rho) \sim \dfrac{\frac{nH}{n-1}}{\sqrt{1 - (\frac{nH}{n-1})^2}} \, \rho.$$
\end{thm}\bigskip

The generating curves are obtained by symmetries from the curves
$(<)$ (standing for $H < \frac{n-1}{n}$) which appear in
Figures~\ref{F-rot-1}-\ref{F-rot-3}.\bigskip

\newpage

\begin{thm}[Rotation $H$-hypersurfaces with $H > \frac{n-1}{n}$]\label{T-h3-rg}
$ $
    \begin{enumerate}
    \item When $d=0$, the hypersurface $\cK_H$ is
        compact and diffeomorphic to an $n$-dimensional sphere. It
        is generated by a compact, simple, strictly convex curve.

    \item When $d>0$, the hypersurface $\cU_H$ is
        complete, embedded and periodic in the $\R$-direction. It looks like
        an unduloid and is contained in a domain of the form $B(0,b)
        \setminus B(0,a) \times \R$, for some constants $0 < a < b$,
        depending on $H$ and $d$.

    \item When $d<0$, the hypersurface $\cN_H$ is
        complete and periodic in the $\R$-direction. It has
        self-intersections, looks like
        a nodoid and is contained in a domain of the form $B(0,b)
        \setminus B(0,a) \times \R$, for some constants $0 < a < b$
        depending on $H$ and $d$.
    \end{enumerate}
\end{thm}\bigskip

The generating curves are obtained by symmetries from the curves
$(>)$ (standing for $H > \frac{n-1}{n}$) which appear in
Figures~\ref{F-rot-1}-\ref{F-rot-3}. \bigskip

\textbf{Remarks}\bigskip

1.~ Constant mean curvature rotation hypersurfaces with $H >
\frac{n-1}{n}$ were obtained in \cite{HH89} and \cite{PR99}.\medskip

2.~ The hypersurfaces $\cS_H$ and the upper (lower) halves of the
cylinders $\cC_H$ in Theorems~\ref{T-h3-r0} and \ref{T-h3-rp} are
stable (as vertical graphs).
\bigskip

\subsection{Proofs of Theorem \ref{T-h3-r0} - \ref{T-h3-rg}}

The proofs follow from an analysis of the asymptotic behaviour of
$I_m(t)$ (Formula (\ref{E-h3-rot2})) when $t$ goes to infinity and
from an analysis of the signs of the functions $nH I_{n-1}(t) + d$,
$M_{H,d}(t)$ and $P_{H,d}(t)$ (Formulas (\ref{E-h3-4})), using the
tables which appear below.
\bigskip

When $d=0$, using (\ref{E-h3-rot1}) one can show that
$\ddot{\lambda} > 0$ and conclude that the generating curve is
strictly convex. When $d \le 0$, the formula for $\ddot{\lambda}$
also shows that the curvature extends continuously at the vertical
points.\bigskip

\noi \textbf{Proof of Theorem \ref{T-h3-r0}}\bigskip

Assume $H = \frac{n-1}{n}$. \bigskip

When $d=0$, the functions $M_{H,0}$ and $P_{H,0}$ are non-negative
and vanish at $t=0$. Near $0$ we have $Q_{H,0}(t) \sim Ht$ and hence
$\lambda_{H,0}(\rho) = \int_0^{\rho} Q_{H,0}(t) \, dt \sim
\frac{H}{2}\rho^2$.
\bigskip

When $d>0$, the function $Q_{H,d}$ exists on an interval $]a_{H,d},
\infty[$ for some constant $a_{H,d} > 0$ and the integral
$\int_{a_{H,d}}^{\rho}Q_{H,d}(t) \, dt$ converges at
$a_{H,d}$.\bigskip

When $d<0$, the function $Q_{H,d}$ exists on an interval
$]\alpha_{H,d}, \infty[$ for some constant $\alpha_{H,d}
> 0$, changes sign from negative to positive, the integral
$\int_{\alpha_{H,d}}^{\rho}Q_{H,d}(t) \, dt$ converges at
$\alpha_{H,d}$ and the curve has a vertical tangent at this point.
The generating curve can be extended by symmetry to a complete curve
with one self-intersection.\bigskip

Using the recurrence relations for the functions $I_m(t)$ one can
determine their asymptotic behaviour at infinity and deduce the
precise exponential growth of the height function $\lambda(\rho)$.

\hfill \qed
\bigskip

\noi \textbf{Proof of Theorem \ref{T-h3-rp}}\bigskip

Assume $0 < H < \frac{n-1}{n}$. \bigskip

When $d=0$, the functions $M_{H,0}$ and $P_{H,0}$ are non-negative
and vanish at $t=0$. Near $0$ we have $Q_{H,0}(t) \sim Ht$ and hence
$\lambda_{H,0}(\rho) = \int_0^{\rho} Q_{H,0}(t) \, dt \sim
\frac{H}{2}\rho^2$. \bigskip

When $d>0$, the function $Q_{H,d}$ exists on an interval $]a_{H,d},
\infty[$ for some constant $a_{H,d} > 0$ and the integral
$\int_{a_{H,d}}^{\rho}Q_{H,d}(t) \, dt$ converges at
$a_{H,d}$.\bigskip

When $d<0$, the function $Q_{H,d}$ changes sign from negative to
positive, exists on an interval $]\alpha_{H,d}, \infty[$ for some
constant $\alpha_{H,d} > 0$, the integral
$\int_{\alpha_{H,d}}^{\rho}Q_{H,d}(t) \, dt$ converges at
$\alpha_{H,d}$ and the generating curve has a vertical tangent at
this point. The generating curve can be extended by symmetry to a
complete curve with one self-intersection.\bigskip

Using the recurrence relations for the functions $I_m(t)$ one can
determine their asymptotic behaviour at infinity and deduce the
precise linear growth of the height function $\lambda(\rho)$.

\hfill \qed
\bigskip

\noi \textbf{Proof of Theorem \ref{T-h3-rg}}\bigskip

Assume $H > \frac{n-1}{n}$. \bigskip

When $d=0$, $Q_{H,0}(t)$ exists on some interval $]0,a_{H,0}[$ for
some positive $a_{H,0}$ and the integral $\lambda_{H,0}(\rho) =
\int_0^{\rho} Q_{H,0}(t) \, dt$ converges at $0$ and at $a_{H,0}$.
The generating curve has a horizontal tangent at $0$ and a vertical
tangent at $a_H$. It can be extended by symmetries to a closed
embedded convex curve.
\bigskip

When $d>0$, the function $Q_{H,d}(t)$ exists on an interval
$]b_{H,d},c_{H,d}[$ for some constants $0 < b_{H,d} < c_{H,d}$ and
the integral converges at the limits of this interval. The
generating curve at these points is vertical. It can be extended by
symmetry to a complete embedded periodic curve (unduloid). \bigskip

When $d<0$, the function $Q_{H,d}(t)$ exists on an interval
$]\beta_{H,d},\gamma_{H,d}[$ for some constants $0 < \beta_{H,d} <
\gamma_{H,d}$, changes sign from negative to positive and the
integral converges at the limits of this interval. The generating
curve at these points is vertical. The generating curve can extended
by symmetries to a complete periodic curve with self-intersections
(nodoid).

\hfill \qed
\bigskip

\textbf{Remark.}~ We note that the integrand $Q_{H,d}(t)$ in
(\ref{E-h3-rot3}) is an increasing function of $H$ for $t$ and $d$
fixed. This fact provides the relative positions of the curves
$\lambda_{H,d}(\rho)$ when $\rho$ and $d$ are fixed. The curve
corresponding to $H > \frac{n-1}{n}$ is above the curve
corresponding to $H = \frac{n-1}{n}$ which is above the curve
corresponding to $H < \frac{n-1}{n}$. See Figures \ref{F-rot-1} to
\ref{F-rot-3}.\bigskip

The above sketches of proof can be completed using the details
below.\bigskip

\noib We have the following relations for the functions $I_m$,

\begin{equation}\label{E-h3-rot2a}
\left\{%
\begin{array}{llll}
  m=0 & I_0(t) & = & t, \\
  m=1 & I_1(t) & = & \cosh(t) - 1, \\
  m=2 & 2I_2(t) & = & \sinh(t) \cosh(t) - t, \\
  m=3 & 3 I_3(t) & = & \sinh^2(t) \cosh(t) - 2(\cosh(t)-1), \\
  m\ge 2 & m I_m(t) & = & \sinh^{m-1}(t) \cosh(t) - (m-1) I_{m-2}(t). \\
\end{array}
\right.
\end{equation}\bigskip

For $m\ge 5$, the asymptotic behavior of $I_m(t)$ near infinity is
given by,

\begin{equation}\label{E-h3-rot2b}
\left\{%
\begin{array}{lll}
  m I_m(t) & = & \sinh^{m-3}(t) \cosh(t) \big(\sinh^2(t) - \frac{m-1}{m-2}\big)
  + O(e^{(m-4)t}),\\
  m I_m(t) & = & \sinh^{m-1}(t) \cosh(t) \big( 1 + O(e^{-2t})\big). \\
\end{array}
\right.
\end{equation}\bigskip

The same holds for $m=4$ with remainder term $O(t)$ in the first
relation. \bigskip

\noib The derivative of $P_{H,d}$ is positive for $t$ positive. The
behaviour of the function $P_{H,d}(t)$ is summarized in the
following table.
\bigskip

\begin{equation*}
    \begin{array}{|c|ccc|}
\hline
n\ge 2   &  & 0 < H &\\
\hline
  t & 0 &  & \infty  \\
\hline
  \partial_t P_{H,d} &  & + &  \\
  \hline
  P_{H,d}(t) & d & \nearrow & \infty  \\
   \hline
\end{array}
\end{equation*}\bigskip

\noib The derivative of $M_{H,d}$ is given by $\partial_t M_{H,d}(t)
= (n-1) \sinh^{n-1}(t) \big( \coth(t) - \frac{nH}{n-1}$\big). For $H
> \frac{n-1}{n}$, we denote by $C_H$ the number such that
$\coth(C_H) = \frac{nH}{n-1}$. The behaviour of the function
$M_{H,d}(t)$ is summarized in the following tables.

\begin{equation*}
    \begin{array}{|c|ccc||ccccc|}
    \hline
 n=2  &  &   0 < H \le \frac{1}{2} &&&& H > \frac{1}{2}&&\\
\hline
  t & 0 &  & \infty & 0&&C_H&& \infty \\
\hline
  \partial_t M_{H,d} &  & + & &&  +&0&-&\\
\hline
  M_{H,d}(t)& -d & \nearrow &
\left\{%
\begin{array}{lr}
  \infty, & H < \frac{1}{2} \\[4pt]
  1-d, & H = \frac{1}{2}\\
\end{array}
\right.
  & -d  & \nearrow & f_H(d) &\searrow & -\infty \\
  \hline
\end{array}
\end{equation*}\bigskip

\begin{equation*}
    \begin{array}{|c|ccc||ccccc|}
    \hline
n\ge 3   &&  0 < H \le \frac{n-1}{n}  && &&H > \frac{n-1}{n}&&\\
\hline
  t & 0 &  & \infty &0&&C_H&& \infty \\
\hline
  \partial_t M_{H,d} &  & + && &+&  0&-&\\
\hline
  M_{H,d}(t) & -d & \nearrow & \infty & -d  & \nearrow & f_H(d) &\searrow & -\infty \\
\hline
\end{array}
\end{equation*}\bigskip

where $f_H(d) := M_{H,d}(C_H) = \sinh^{n-1}(C_H) - nH I_{n-1}(C_H)
-d$. \bigskip

The signs and zeroes of the functions $M_{H,d}(t)$ and $P_{H,d}(t)$
when $d \not = 0$ are summarized in the following charts, together
with the existence domain of the function $Q_{H,d}$.\bigskip

When $d > 0$, we have

\begin{equation*}
    \begin{array}{|c|cccccc|}
\hline n=2 &&
\begin{array}{c}
0 < H < \frac{1}{2},\\
H=\frac{1}{2},
\end{array}
&&
\begin{array}{c}
0 < d \\
0<d<1,
\end{array}
&&\\
\hline n \ge 3 &&&
\begin{array}{c}
0 < H \le \frac{n-1}{n}\\
0 < d
\end{array}
&&&\\
\hline
t & 0 & & a_{H,d} && & \infty \\
\hline
M_{H,d} && - & 0 & +& & \\
\hline
P_{H,d} && + & + & +& & \\
\hline
Q_{H,d} &&  \not \exists & +\infty &\exists && \\
\hline
    \end{array}
\end{equation*}\bigskip

\begin{equation*}
    \begin{array}{|c|ccccccc|}
\hline n \ge 2 &&&&
\begin{array}{c}
H > \frac{n-1}{n}\\
0 < d < D_H
\end{array}
&&&\\
\hline
t & 0 & & b_{H,d} & C_{H} & c_{H,d} & & \infty \\
\hline
M_{H,d}& &-& 0 & + & 0& -& \\
\hline
P_{H,d}& &+& & + & & +& \\
\hline
Q_{H,d}& &\not \exists & +\infty & \exists & +\infty & \not \exists & \\
\hline
    \end{array}
\end{equation*}\bigskip

where $D_H := \sinh^{n-1}(C_H) - nH I_{n-1}(C_H)$. \bigskip


\newpage

When $d<0$, we have the following tables.\bigskip

\begin{equation*}
    \begin{array}{|c|ccccccc|}
\hline n \ge 2 &&&&
\begin{array}{c}
0 < H \le \frac{n-1}{n}\\
d < 0
\end{array}
&&&\\
\hline
t & 0 & &  & \alpha_{H,d} & & & \infty \\
\hline
M_{H,d}& && + &  & +& & \\
\hline
P_{H,d}& && - & 0 & +& & \\
\hline
Q_{H,d}& & & \not \exists & - \infty & \exists &  & \\
\hline
    \end{array}
\end{equation*}\bigskip

Note that the function $Q_{H,d}$ changes sign from negative to
positive when $t$ goes from $\alpha_{H,d}$ to infinity.\bigskip

\begin{equation*}
    \begin{array}{|c|ccccccc|}
\hline n \ge 2 &&&&
\begin{array}{c}
H > \frac{n-1}{n}\\
d < 0
\end{array}
&&&\\
\hline
t & 0 & & \gamma_{H,d} &  & \beta_{H,d} & & \infty \\
\hline
M_{H,d}& &+& + & +& 0& -& \\
\hline
P_{H,d}& &-& 0 & + & +& +& \\
\hline
Q_{H,d}& &\not \exists & -\infty & \exists & +\infty & \not \exists & \\
\hline
    \end{array}
\end{equation*}\bigskip

Note that the function $Q_{H,d}$ changes sign from negative to
positive when $t$ goes from $\gamma_{H,d}$ to $\beta_{H,d}$.\bigskip


\subsection{Translation invariant $H$-hypersurfaces
in $\HH^n \times \R$}\label{SS-dim3-Htransl}
\bigskip

\subsubsection{Translation hypersurfaces}
\bigskip

\noib \textbf{Definitions and Notations.}~ We consider $\gamma$ a
geodesic through $0$ in $\HH^n$ and the totally geodesic vertical
plane $\V = \gamma \times \R = \ens{(\gamma (\rho),t)}{(\rho ,t) \in
\R \times \R}$ where $\rho$ is the signed hyperbolic distance to $0$
on $\gamma$.
\bigskip

Take $\PP$ a totally geodesic hyperplane in $\HH^n$, orthogonal to
$\gamma$ at $0$. We consider the hyperbolic translations with
respect to the geodesics $\delta$ through $0$ in $\PP$. We shall
refer to these translations as translations with respect to $\PP$.
These isometries of $\HH^n$ extend ``slice-wise'' to isometries of
$\HH^n \times \R$. \bigskip

In the vertical plane $\V$, we consider the curve $c(\rho) := \big(
\tanh(\rho /2), \mu(\rho)\big)$. \bigskip

In $\HH^n \times \{\mu(\rho)\}$, we translate the point $c(\rho)$ by
the translations with respect to $\PP\times \{\mu(\rho)\}$ and we
get the equidistant hypersurface $\PP_{\rho}$ passing through
$c(\rho)$, at distance $\rho$ from $\PP\times \{\mu(\rho)\}$. The
curve $c$ then generates a \emph{translation hypersurface} $M =
\cup_{\rho}\PP_{\rho}$ in $\HH^n \times \R$.\bigskip

\noib \textbf{Principal curvatures.}~ The principal directions of
curvature of $M$ are the tangent to the curve $c$ in $\V$ and the
directions tangent to $\PP_{\rho}$. The corresponding principal
curvatures with respect to the unit normal pointing upwards are
given by

\begin{equation*}\label{E-tra-1}
\left\{%
\begin{array}{lll}
    k_{\V} & = & \ddot{\mu}(\rho) \big( 1 + \dot{\mu}^2(\rho) \big)^{-3/2}, \\
    k_{\PP} & = &  \dot{\mu}(\rho) \big( 1 + \dot{\mu}^2(\rho) \big)^{-1/2}
    \tanh(\rho). \\
\end{array}%
\right.
\end{equation*}\bigskip

The first equality comes from the fact that $\V$ is totally geodesic
and flat. The second equality follows from the fact that
$\PP_{\rho}$ is totally umbilic and at distance $\rho$ from
$\PP\times \{\mu(\rho)\}$ in $\HH^n \times \{\mu(\rho)\}$. \bigskip

\noib \textbf{Mean curvature.}~ The mean curvature of the
translation hypersurface $M$ associated with $\mu$ is given by



\begin{equation}\label{E-tra-3}
n H(\rho) \cosh^{n-1}(\rho) = \partial_{\rho} \Big(
\cosh^{n-1}(\rho) \dot{\mu}(\rho) \big( 1 + \dot{\mu}^2(\rho)
\big)^{-1/2} \Big).
\end{equation}\bigskip

\subsubsection{Constant mean curvature translation hypersurfaces}
\bigskip


We may assume that $H \ge 0$. The generating curves of translation
hypersurfaces with constant mean curvature $H$ are given by the
differential equation

\begin{equation}\label{E-tra-4}
\dot{\mu}(\rho) \big( 1 + \dot{\mu}^2(\rho) \big)^{-1/2}
\cosh^{n-1}(\rho) = n H \int_0^{\rho} \cosh^{n-1}(t) \, dt + d
\end{equation}

for some integration constant $d$.\bigskip

Minimal translation hypersurfaces have been studied in \cite{Sa08,
ST08} in dimension $2$ and in \cite{BS08a} in higher dimensions.
Constant mean curvature ($H \not = 0$) translation hypersurfaces
have been treated in \cite{Sa08} in dimension $2$. The purpose of
the present section is to investigate the higher dimensional
translation $H$-hypersurfaces.\bigskip

\textbf{Notations.}~ For later purposes, we introduce some
notations.\bigskip

\noib For $m \ge 0$, we define the functions

\begin{equation}\label{E-tra-5a}
J_m(r) := \int_0^r \cosh^m(t) \, dt.
\end{equation}\bigskip

\noib For $H > 0$ and $d \in \R$, we introduce the functions,

\begin{equation}\label{E-tra-11}
\left\{%
\begin{array}{lll}
R_{H,d}(t) & = & \cosh^{n-1}(t) - nH J_{n-1}(t) -d ,\\
S_{H,d}(t) & = & \cosh^{n-1}(t) + nH J_{n-1}(t) +d ,\\
T_{H,d}(t) & = & \big[ nH J_{n-1}(t) + d \big]
\big[R_{H,d(t)} S_{H,d}(t) \big]^{-1/2}.\\
\end{array}%
\right.
\end{equation}\bigskip

We note from (\ref{E-tra-4}) that $\dot{\mu}(t)$ has the sign of
$nHJ_{n-1}(t) +d$. It follows that $\mu$ is given (up to an additive
contant) by

\begin{equation*}\label{E-tra-8a}
\mu_{H,d} (\rho) = \int_{\rho_0}^{\rho} \big[ nH J_{n-1}(t) + d
\big] \big[ \cosh^{2n-2}(t) - \big( nH J_{n-1}(t) + d \big)^2
\big]^{-1/2} \, dt
\end{equation*}

or, using the above notations,

\begin{equation}\label{E-tra-8}
\mu_{H,d} (\rho) = \int_{\rho_0}^{\rho} \big[ nH J_{n-1}(t) + d
\big] \big[R_{H,d(t)} \, S_{H,d}(t) \big]^{-1/2} \, dt =
\int_{\rho_0}^{\rho} T_{H,d}(t) \, dt \, ,
\end{equation}

where the integration interval $[\rho_0, \rho]$ is contained in the
interval in which the square root exists. The existence and
behaviour of the function $\mu_{H,d}$ depend on the signs of the
functions $nH J_{n-1}(t) + d$, $R_{H,d}(t)$ and
$S_{H,d}(t)$.\bigskip

For $H=\frac{n-1}{n}$, we give a complete description of the
corresponding translation $H$-hypersurfaces. For $0 < H <
\frac{n-1}{n}$, we prove the existence of a complete non-entire
$H$-graph with infinite boundary data and infinite asymptotic
behaviour. The other cases can be treated similarly using the tables
below. \bigskip

\begin{thm}[Translation $H$-hypersurfaces, with $n\ge 3$ and $H =
\frac{n-1}{n}$]\label{T-tra-1}
$ $
\begin{enumerate}
    \item When $d=0$, $\cT_0$ is a complete embedded smooth
    hypersurface generated by a compact, simple, strictly convex curve.
    The hypersurface is symmetric with respect
    to a horizontal hyperplane and the parts above and below this hyperplane
    are vertical graphs. The hypersurface also admits a vertical
    symmetry. The asymptotic boundary of $\cT_0$ is topologically
    a cylinder.

    \item When $0 < d < 1$, the hypersurface $\cT_d$ is similar to
    $\cT_0$ except that it is not smooth.

    \item When $d \le -1$, $\cT_d$ is a smooth complete immersed
    hypersurface with self-intersections and horizontal symmetries.
    The asymptotic boundary of $\cT_d$ is topologically a cylinder.

    \item When $-1 < d < 0$, the hypersurface $\cT_d$ looks like
    $\cT_{-1}$ except that it is not smooth.
\end{enumerate}
\end{thm}

\begin{pb1-figs}
\begin{figure}[h]
\begin{center}
\begin{minipage}[c]{6.5cm}
    \includegraphics[width=6.5cm]{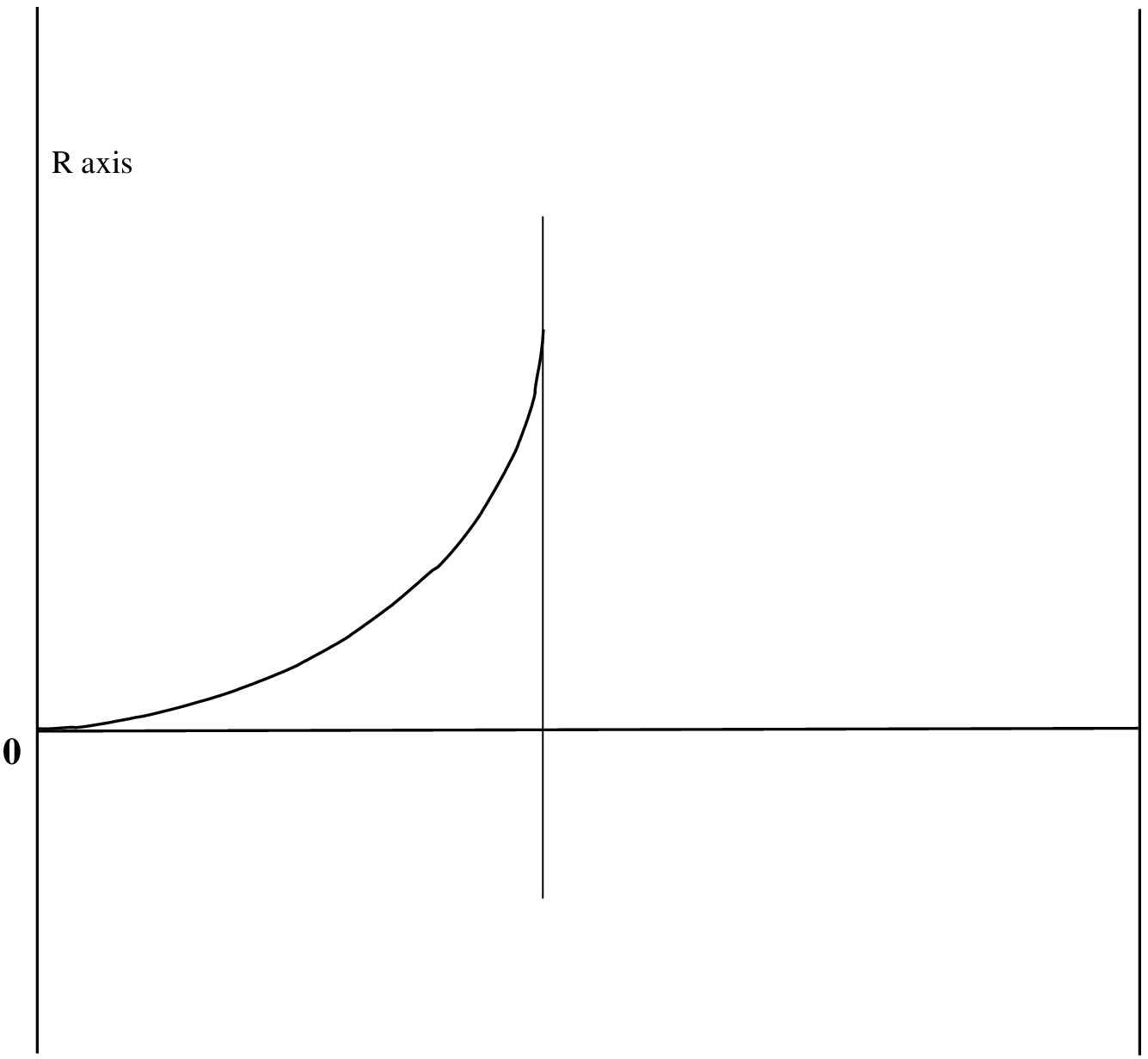}
    \caption[$n \ge 3, H =\frac{n-1}{n}, d = 0$]{$n \ge 3, H =\frac{n-1}{n}, d = 0$}
    \label{F-tra-3a}
\end{minipage}\hfill
\begin{minipage}[c]{6.5cm}
    \includegraphics[width=6.5cm]{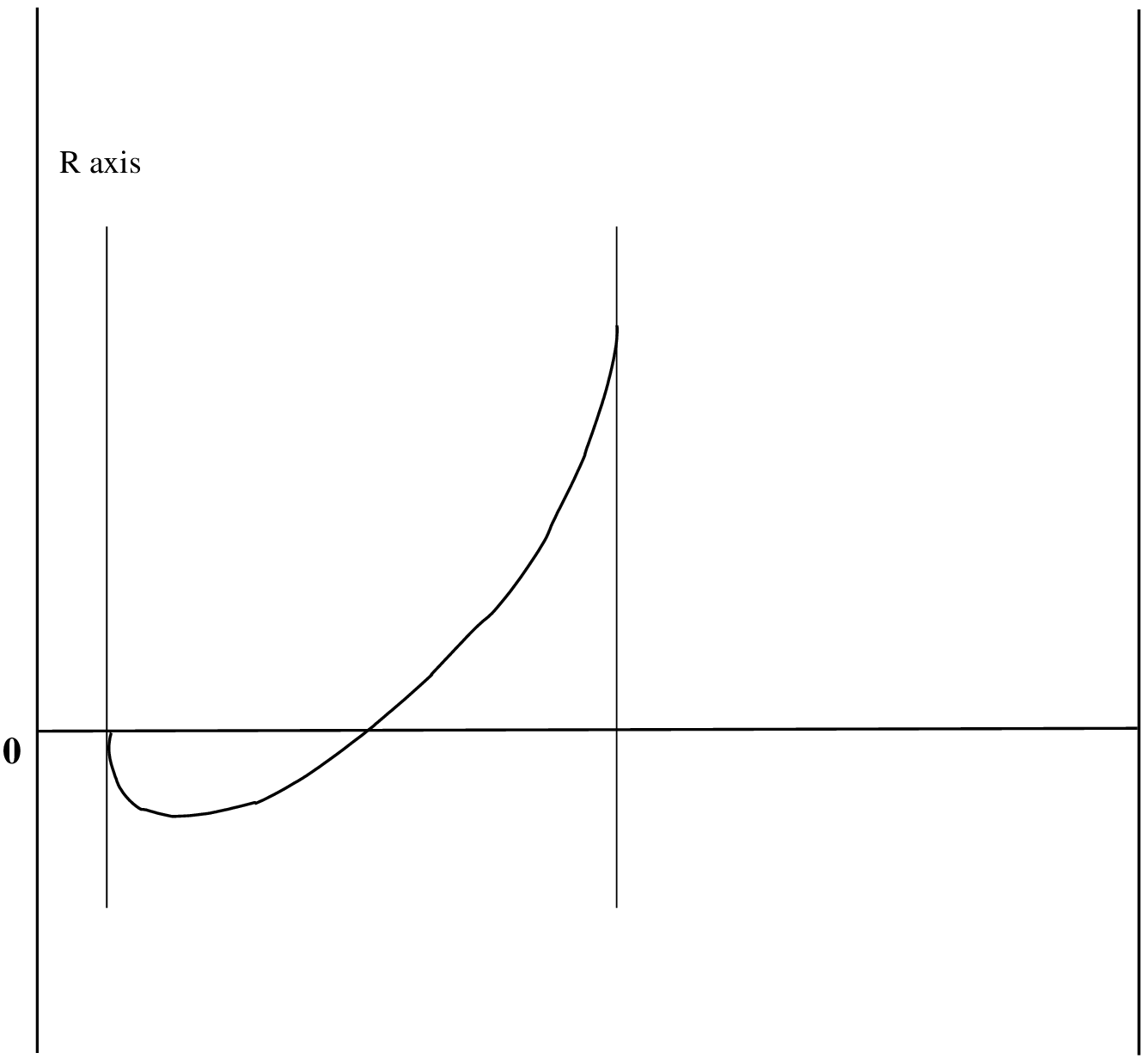}
    \caption[$n \ge 3, H =\frac{n-1}{n}, d < -1$]{$n \ge 3, H =\frac{n-1}{n}, d < -1$}
    \label{F-tra-3c}
\end{minipage}\hfill
\end{center}
\end{figure}
\end{pb1-figs}\bigskip

\textbf{Remark.}~ When $d \ge 1$, the differential equation
(\ref{E-tra-4}) does not have solutions.

\begin{pb1-figs}
\begin{figure}[h]
\begin{center}
\begin{minipage}[c]{6.5cm}
    \includegraphics[width=6.5cm]{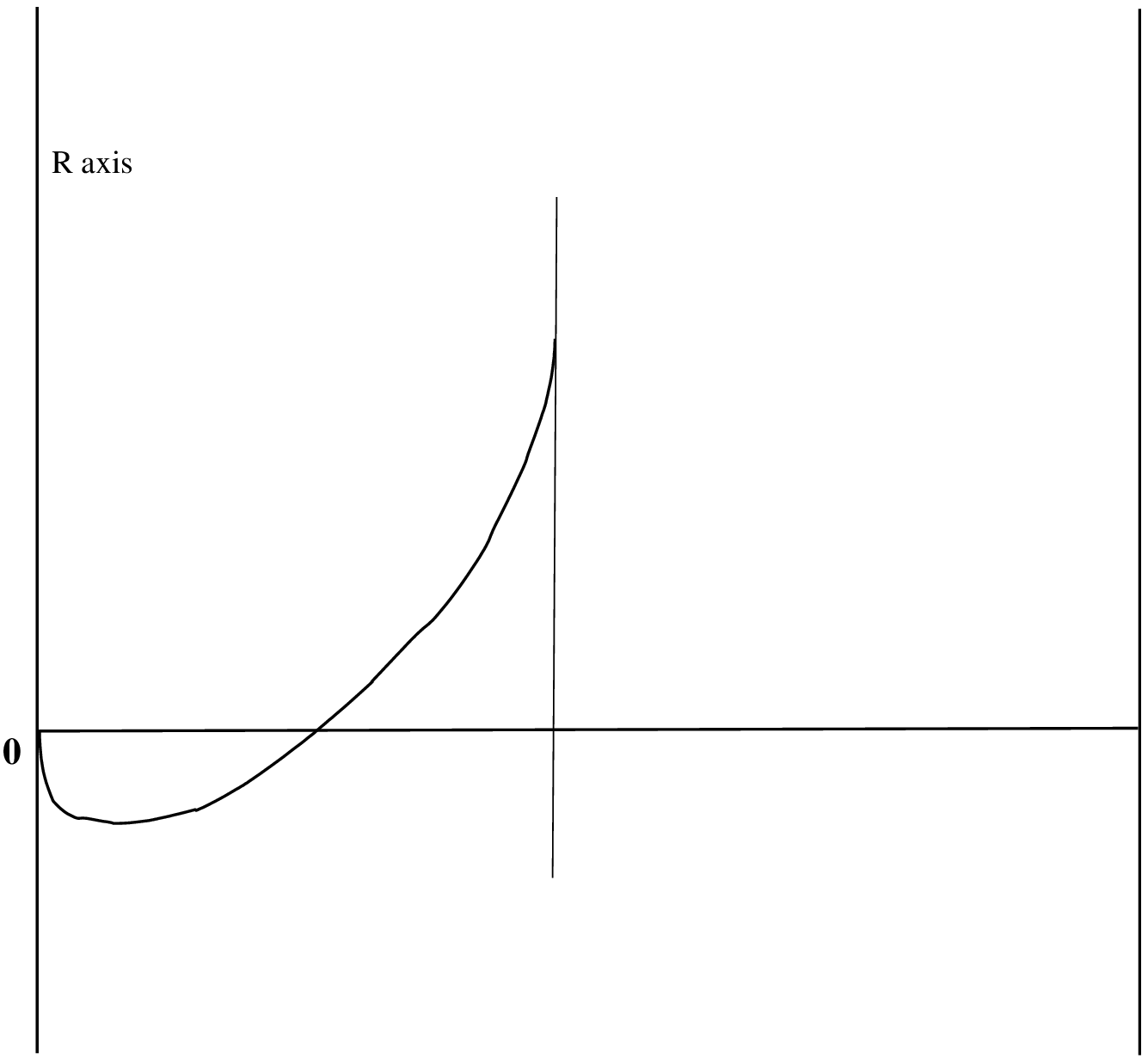}
    \caption[$n \ge 3, H =\frac{n-1}{n}, d = - 1$]{$n \ge 3, H =\frac{n-1}{n}, d = - 1$}
    \label{F-tra-3d}
\end{minipage}\hfill
\begin{minipage}[c]{6.5cm}
    \includegraphics[width=6.5cm]{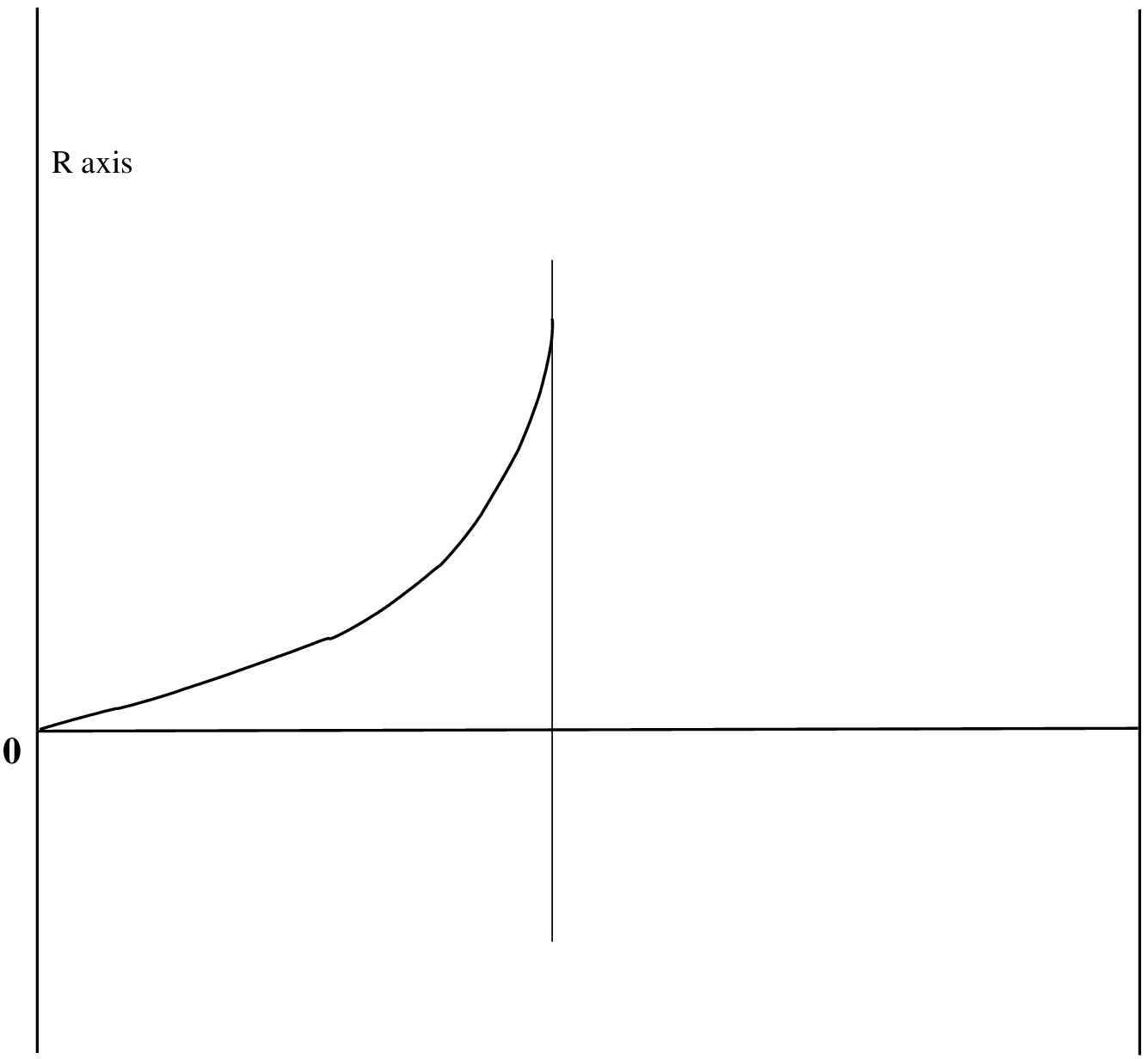}
    \caption[$n \ge 3, H =\frac{n-1}{n}, 0 < d < 1$]{$n \ge 3,
    H =\frac{n-1}{n}, 0 < d < 1$}
    \label{F-tra-3b}
\end{minipage}\hfill
\end{center}
\end{figure}
\end{pb1-figs}

\begin{pb1-figs}
\begin{figure}[ht]
\begin{center}
\begin{minipage}[c]{6.5cm}
\includegraphics[width=6.5cm]{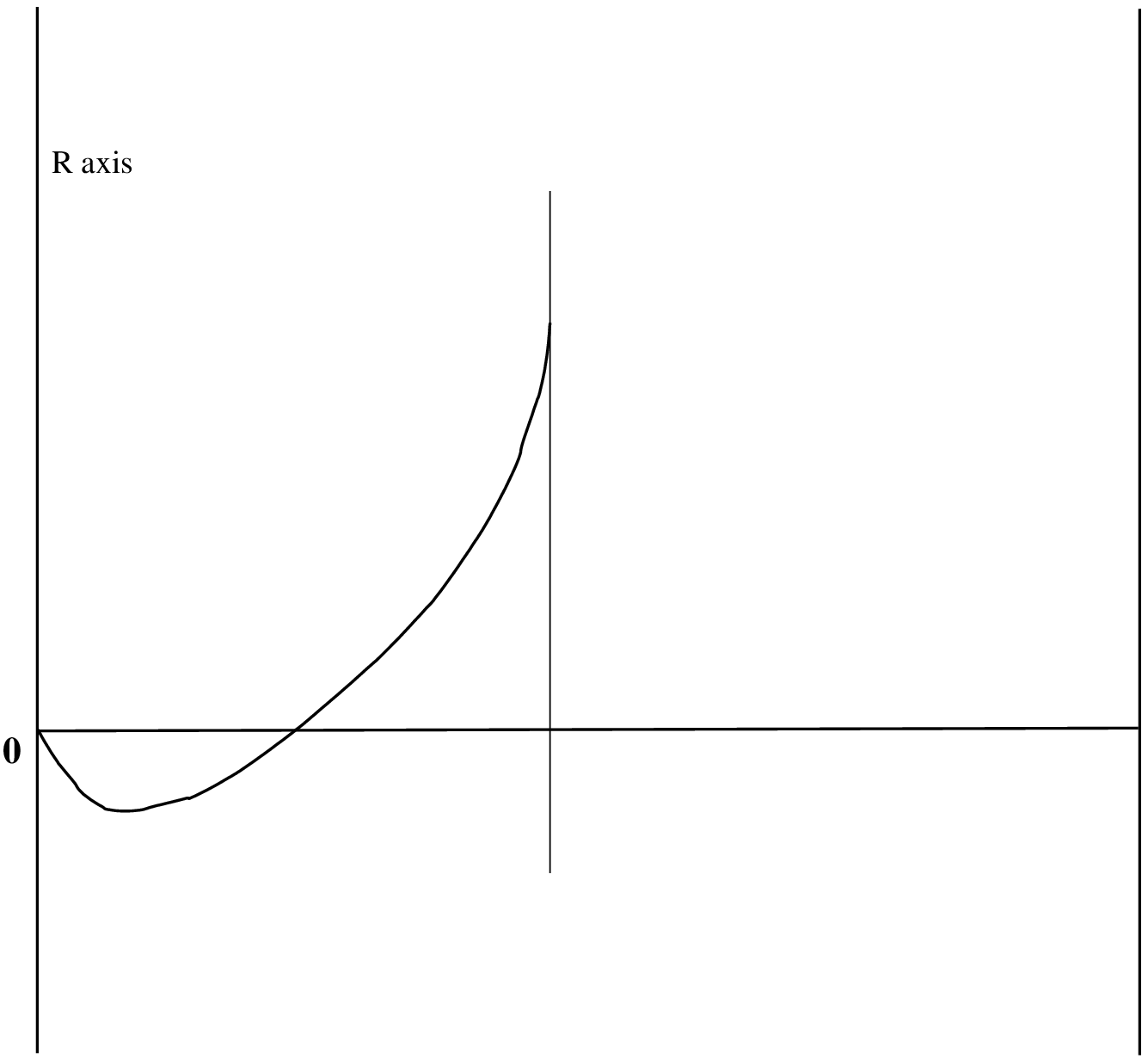}
    \caption[$n \ge 3, H =\frac{n-1}{n}, -1 < d < 0$]{$n \ge 3,
    H =\frac{n-1}{n}, -1 < d < 0$}
    \label{F-tra-3e}
\end{minipage}\hfill
\begin{minipage}[c]{6.5cm}
\includegraphics[width=6.5cm]{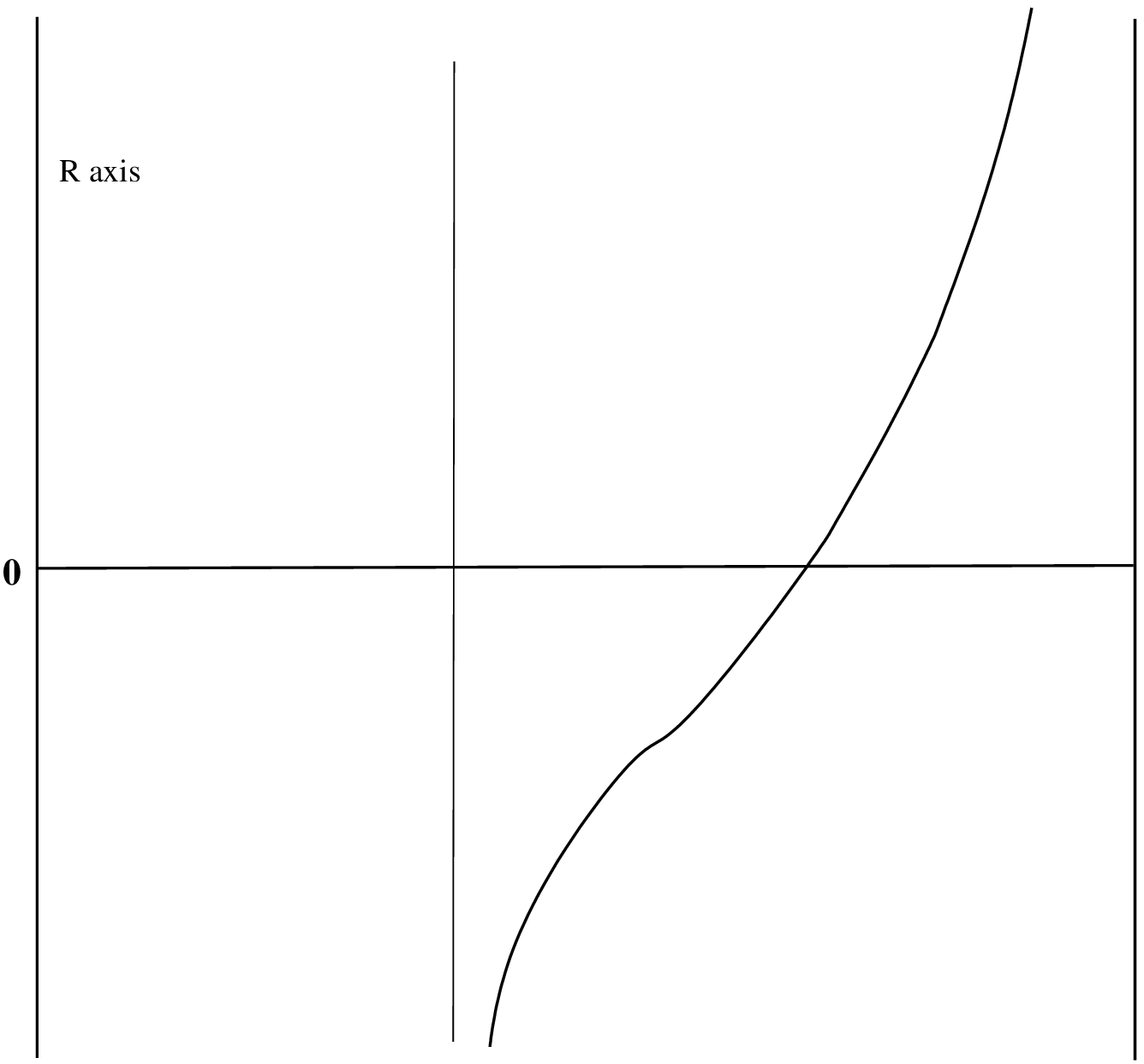}
    \caption[$n \ge 2, H <\frac{n-1}{n}$]{$n \ge 2,
    H < \frac{n-1}{n}$}
    \label{F-tra-33}
\end{minipage}\hfill
\end{center}
\end{figure}\bigskip
\end{pb1-figs}

\begin{thm}[Complete $H$-graph with infinite boundary data]\label{T-tra-2}
$ $\\[2pt]
There exists a complete translation hypersurface $\cT_{H}$, with $0
< H < \frac{n-1}{n}$, such that
\begin{enumerate}
\item $\cT_H$ is a complete monotone vertical $H$-graph over the non mean
convex side of an equidistant hypersurface $\Gamma\subset \HH^n$
with mean curvature $\frac{nH}{n-1},$
\item $\cT_H$ takes infinite boundary value data on $\Gamma$ and infinite
asymptotic boundary data.
\end{enumerate}
\end{thm}\bigskip

\subsection{Proof of Theorem \ref{T-tra-1}}

\bigskip

The proof of Theorem \ref{T-tra-1} follows from an analysis of the
asymptotic behaviour of the functions $J_m(t)$ (Formula
(\ref{E-tra-5a})) when $t$ goes to infinity and from an analysis of
the signs of the functions $R_{H,d}$ and $S_{H,d}$ (Formulas
(\ref{E-tra-11})) depending on the signs of $H - \frac{n-1}{n}$ and
$d$.\bigskip

\noib We have the relations

\begin{equation}\label{E-tra-10}
\left\{%
\begin{array}{lll}
J_0(t) & = & t ,\\
J_1(t) & = & \sinh(t) ,\\
2J_2(t) & = & \sinh(t) \cosh(t) + t ,\\
3J_3(t) & = & \sinh(t) \cosh^2(t) + 2 J_1(t) ,\\
mJ_m(t) & = & \sinh(t) \cosh^{m-1}(t) + (m-1) J_{m-2}(t), \text{ ~for~ } m\ge 3.\\
\end{array}%
\right.
\end{equation}\bigskip

These relations give us the asymptotic behaviour of the functions
$J_m(t)$ when $t$ tends to infinity. In particular,

$$
m J_m(t) = \sinh(t) \cosh^{m-1}(t) + \frac{m-1}{m-2} \sinh(t)
\cosh^{m-3}(t) + O(e^{(m-4)t}), \text{ ~for~ } m\ge 5
$$

with the remainder term replaced by $O(t)$ when $m=4$.\bigskip

\noib \textbf{The function $S_{H,d}(t)$}\bigskip

For all $H > 0$, the function $S_{H,d}$ increases from $1+d$ to $+
\infty$. Its behaviour is summarized in the following table.

\begin{equation}\label{E-tra-12s}
\begin{array}{|c|ccccc|}
  \hline
   & \text{Case} & 0<H &  & d \ge -1 &  \\
  \hline
  t & 0 &  &  &  & + \infty \\
  \hline
  S_{H,d}(t) & 1+d \ge 0 &  & \nearrow &  & +\infty \\
  \hline \hline
   & \text{Case} & 0<H &  & d < -1 &  \\
  \hline
  t & 0 &  & \alpha_{H,d} &   & + \infty \\
  \hline
  S_{H,d}(t) & 1+d<0 & \nearrow & 0 & \nearrow & + \infty \\
  \hline
\end{array}
\end{equation}\bigskip

\noib \textbf{The function $R_{H,d}(t)$}\bigskip

The derivative of $R_{H,d}(t)$ is given by $\partial_t R_{H,d}(t) =
(n-1) \cosh^{n-1}(t) [\tanh(t) - \frac{nH}{n-1}]$. For $0 < H <
\frac{n-1}{n}$, let $t_H$ be the value such that $\tanh(t_H) =
\frac{nH}{n-1}$.\bigskip

\noi When $H \not = \frac{n-1}{n}$,

\begin{equation}\label{E-tra-asymp}
R_{H,d}(t) \sim \frac{1}{2}(1-\frac{nH}{n-1}) \cosh^{n-2}(t) e^t
\text{ ~near~ } t = + \infty.
\end{equation}\bigskip

\noi When $H = \frac{n-1}{n}$ and when $t$ tends to $+ \infty$,
$R_{H,d}(t)$ tends to $- \infty$ for $n \ge 3$ and to $-d$ for
$n=2$.\bigskip

The behaviour of the function $R_{H,d}(t)$ is summarized in the
following table. \bigskip

\begin{equation}\label{E-tra-12r}
\begin{array}{|c|ccccc|}
  \hline
   & \text{Case} &  & 0<H<\frac{n-1}{n} &  &  \\
  \hline
  t & 0 &  & t_H &  & + \infty \\
  \hline
  R_{H,d}(t) & 1-d & \searrow & R_{H,d}(t_H) & \nearrow & +\infty \\
  \hline
\hline
   & \text{Case} &  & H = \frac{n-1}{n} &  &  \\
  \hline
  t & 0 &  &  &  & + \infty \\
  \hline
  R_{H,d}(t) & 1-d &  & \searrow &  & \left\{%
                                      \begin{array}{cc}
                                     - \infty, & n \ge 3 \\
                                        -d, & n=2 \\
                                      \end{array}
                                        \right. \\
  \hline
  \hline
   & \text{Case} &  & H > \frac{n-1}{n} &  &  \\
  \hline
  t & 0 &  &  &  & + \infty \\
  \hline
  R_{H,d}(t) & 1-d &  & \searrow &  & -\infty \\
  \hline
\end{array}
\end{equation}\bigskip

\textbf{Proof of Theorem \ref{T-tra-1}, continued}\bigskip

We now investigate the behaviour of the solution $\mu$ to Equation
(\ref{E-tra-4}) when $n\ge 3$ and $H=\frac{n-1}{n}$ (for $n=2$, see
\cite{Sa08}).\bigskip

\noi According to Table (\ref{E-tra-12s}), the function $S_{H,d}$
increases from $1+d$ to $+ \infty$ and we have to consider two
cases, \emph{(i)} $d \ge - 1$, in which case $S_{H,d}$ is always
non-negative and \emph{(ii)} $d < -1$, in which case $S_{H,d}$ has
one zero $\alpha_{H,d}$ such that $$\cosh^{n-1}(\alpha_{H,d}) + nH
J_{n-1}(\alpha_{H,d}) + d = 0.$$

\noi According to Table (\ref{E-tra-12r}), the function $R_{H,d}$
decreases from $1-d$ to
$\left\{%
\begin{array}{cc}
- \infty, & n \ge 3 \\
-d, & n=2 \\
\end{array} \right. $, depending on the value of $n$. It follows
that we have two cases, \emph{(i)} $d \ge 1$, in which case the
function $R_{H,d}$ is always non-positive and \emph{(ii)} $d < 1$,
in which case it has one zero $c_{H,d}$ for $n \ge 3$. When it
exists, the zero $c_{H,d}$ satisfies  $$\cosh^{n-1}(c_{H,d}) - nH
J_{n-1}(c_{H,d}) - d = 0.$$
\bigskip

Looking at the equations defining $\alpha_{H,d}$ and $c_{H,d}$ we
see that $\alpha_{H,d} < c_{H,d}$ when they both exist.\bigskip

The behaviour of the function $\mu$ is described in the following
tables, see also Figures~\ref{F-tra-3a} to \ref{F-tra-3e}.\bigskip

\begin{equation}\label{E-tra-15-1}
\begin{array}{|c|ccccccc|}
\hline
  \textbf{Case 1} &  & H=\frac{n-1}{n} &   & d < -1 &  & n\ge 3 &  \\
\hline
  t & 0 &  & \alpha_{H,d} &  & c_{H,d} &  & + \infty \\
\hline
  R_{H,d} &  & + & + & + & 0 & - &  \\
\hline
  S_{H,d} &  & - & 0 & + & + & + &  \\
\hline
  T_{H,d} &  & \not \exists & - \infty & \exists  & + \infty
  & \not \exists &  \\
\hline
\end{array}
\end{equation}\bigskip

The function $\mu$ is given by

$$\mu(\rho) = \int_{\rho_0}^{\rho} T_{H,d}(t) \, dt$$

for $\rho_0, \rho \in [\alpha_{H,d}, c_{H,d}]$ and the integral
exists at both limits. Note that the integrand is negative near the
lower limit while it is positive near the upper limit. \bigskip

When $d=0$, using (\ref{E-tra-4}) one can show that $\ddot{\mu}
> 0$ and conclude that the generating curve is strictly convex. The
formula for $\ddot{\mu}$ also shows that the curvature extends
continuously at the vertical points.\bigskip

The generating curve can be extended by symmetry and periodicity to
give rise to a complete immersed hypersurface with
self-intersections.\bigskip

\begin{equation}\label{E-tra-15-2}
\begin{array}{|c|ccccccc|}
\hline
  \textbf{Case 2} &  & H=\frac{n-1}{n} &   & -1 \le d < 1 &  & n\ge 3 &  \\
\hline
  t & 0 &  &  &  & c_{H,d} &  & + \infty \\
\hline
  R_{H,d} &  &  & + &  & 0 & - &  \\
\hline
  S_{H,d} &  &  & + &  & + & + &  \\
\hline
  T_{H,d} &  &  & \exists &   & + \infty
  & \not \exists &  \\
\hline
\end{array}
\end{equation}\bigskip

The function $\mu$ is given by

$$\mu(\rho) = \int_{0}^{\rho} T_{H,d}(t) \, dt$$

for $\rho_0, \rho \in [0, c_{H,d}]$ and the integral exists at both
ends. Note that the integrand has the sign of $d$ near $0$, with
$\dot{\mu}(0) = d/\sqrt{1-d^2}$ ; it is positive near the upper
bound with $\dot{\mu}(c_{H,d}) = + \infty$. \bigskip

When $d=-1$, the original curve has a vertical tangent at $0$. It
can be extended by symmetry and periodicity to give rise to a
complete immersed hypersurface with self-intersections. \bigskip

When $d=0$, the curve has a horizontal tangent and is strictly
convex (use (\ref{E-tra-4})). It can be extended by symmetry as a
topological circle and gives rise to a complete embedded surface.
\bigskip

When $d \ge 1$, Equation (\ref{E-tra-4}) has no solution.\bigskip

\subsection{Proof of Theorem \ref{T-tra-2}}

Given $n$ and $H$, such that $0<H<\frac{n-1}{n}$, consider the
function $R_{H,d}(t)$ and choose $d_H$ such that $R_{H,d_H}(t_H) =
0$, where $t_H$ is defined by $\tanh (t_H) = \frac{nH}{n-1}$, \ie
$d_H := \cosh^{n-1}(t_H) - nH J_{n-1}(t_H)$.\bigskip

It follows that  $R_{H,d_H}(t) > 0$ for $t > t_H$ and hence the
quantity $nH J_{n-1}(t) +d_H$ does not change sign for $t>t_H$ and
the same is true  for $T_{H,d_H}(t)$.\bigskip

Taking (\ref{E-tra-11}) into account, we choose $\rho_0 > t_H$ and
define the generating curve by Formula (\ref{E-tra-8}).\bigskip


We conclude that $\mu(\rho)$ is well-defined and strictly increasing
for $\rho>t_H$. Moreover, $\mu(\rho)$ goes to $-\infty$, if
$\rho\rightarrow t_H^+.$  Notice that the mean curvature of the
equidistant hypersurface at distance $t_H$ to $\PP$ is $\tanh
(t_H)=\frac{n H}{n-1}$, by the choice of $t_H$.\bigskip

 Now recall that if $ 0<H<\frac{n-1}{n},$ then $R_{H,d}(t)\sim
 \frac{1}{2}(1-\frac{nH}{n-1})\cosh^{n-2}(t) e^t,$ as
 $t\rightarrow \infty.$ From this it follows that $T_{H,d}(t)=O(1)$,
 as $t\rightarrow \infty.$ Thus $\mu(\rho)\rightarrow
 +\infty$, if $\rho\rightarrow \infty$.
 \hfill  \qed

\vspace{1.5cm}

\textbf{Remark.}~ The situation when $n=2$ is similar although the
generating curves are defined on infinite intervals (see Figures
\ref{F-tra-2a} to \ref{F-tra-2e}). The corresponding surfaces have
height functions tending to infinity when $\rho$ tends to infinity.
In particular, the surface $\cT_0$ is a complete smooth entire graph
above $\HH^2$.

\vspace{1cm}


\begin{pb1-figs}
\begin{figure}[h]
\begin{center}
\begin{minipage}[c]{4.5cm}
    \includegraphics[width=4.5cm]{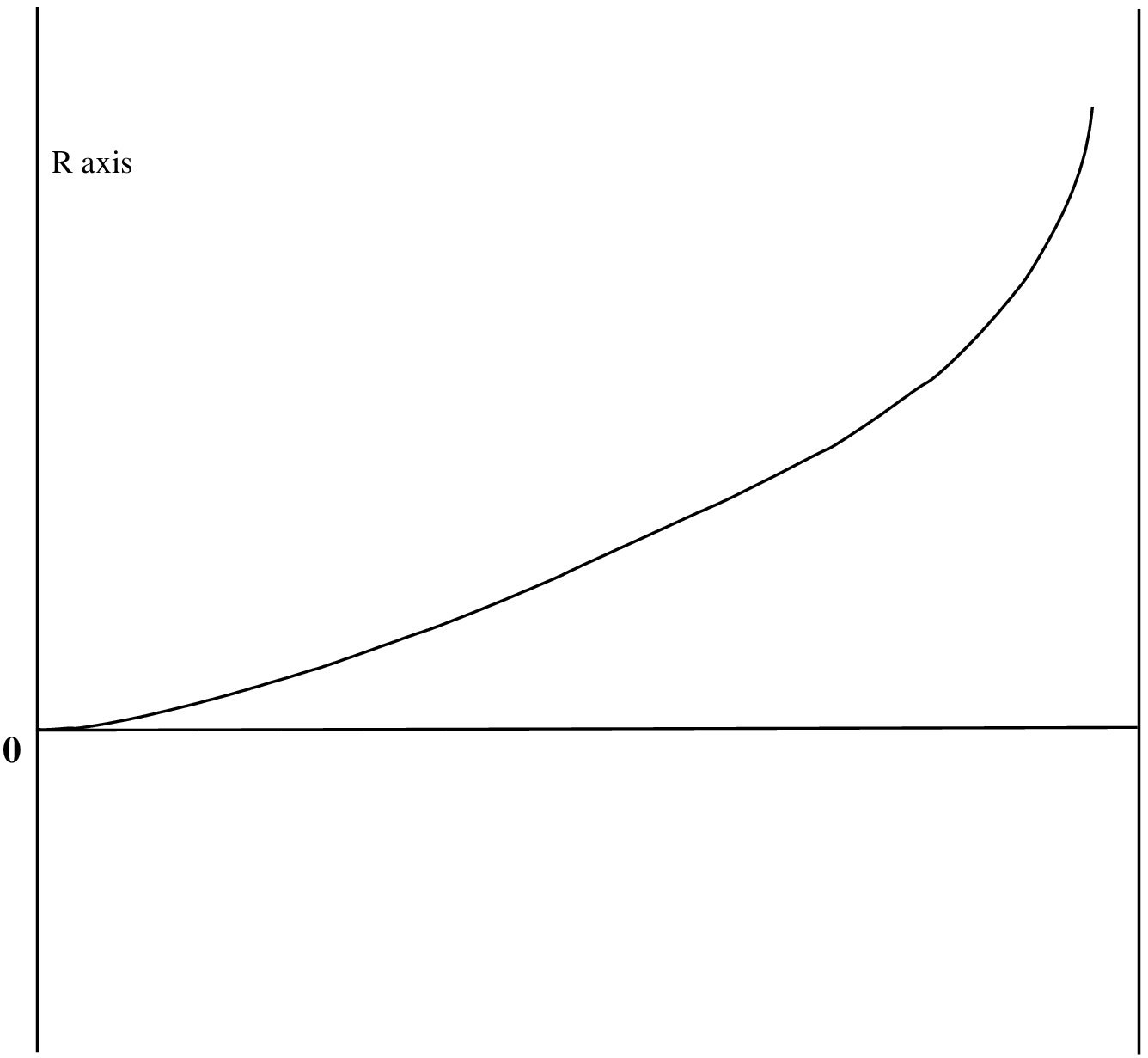}
    \caption[$n=2, H =\frac{1}{2}$, $d = 0$]{$n=2, H =\frac{1}{2}$, \\[3pt] $d = 0$}
    \label{F-tra-2a}
\end{minipage}\hfill
\begin{minipage}[c]{4.5cm}
\includegraphics[width=4.5cm]{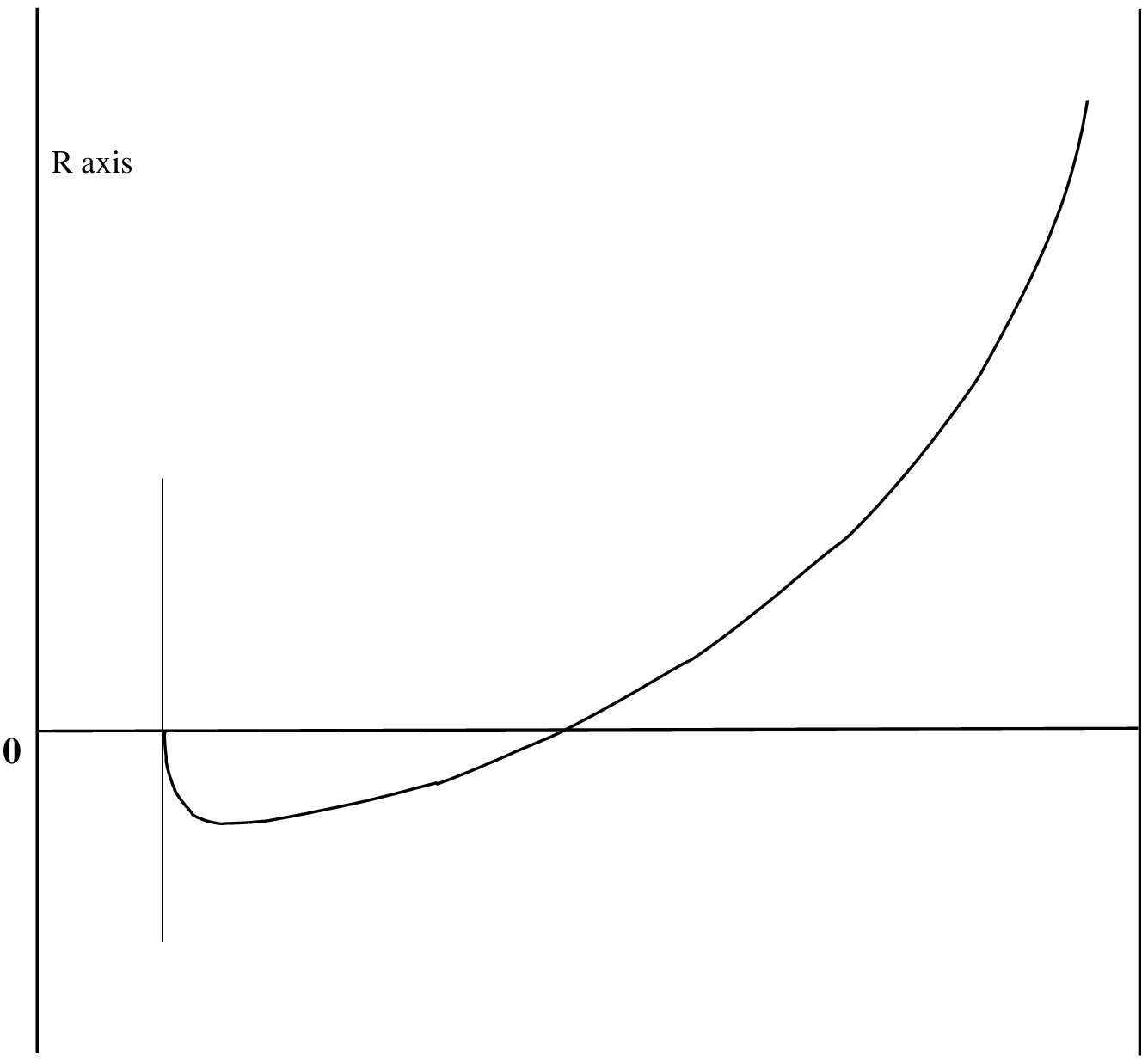}
    \caption[$n=2, H =\frac{1}{2}$, $d < -1$]{$n=2, H =\frac{1}{2}$, \\ $d < -1$}
    \label{F-tra-2c}
\end{minipage}\hfill
\begin{minipage}[c]{4.5cm}
    \includegraphics[width=4.5cm]{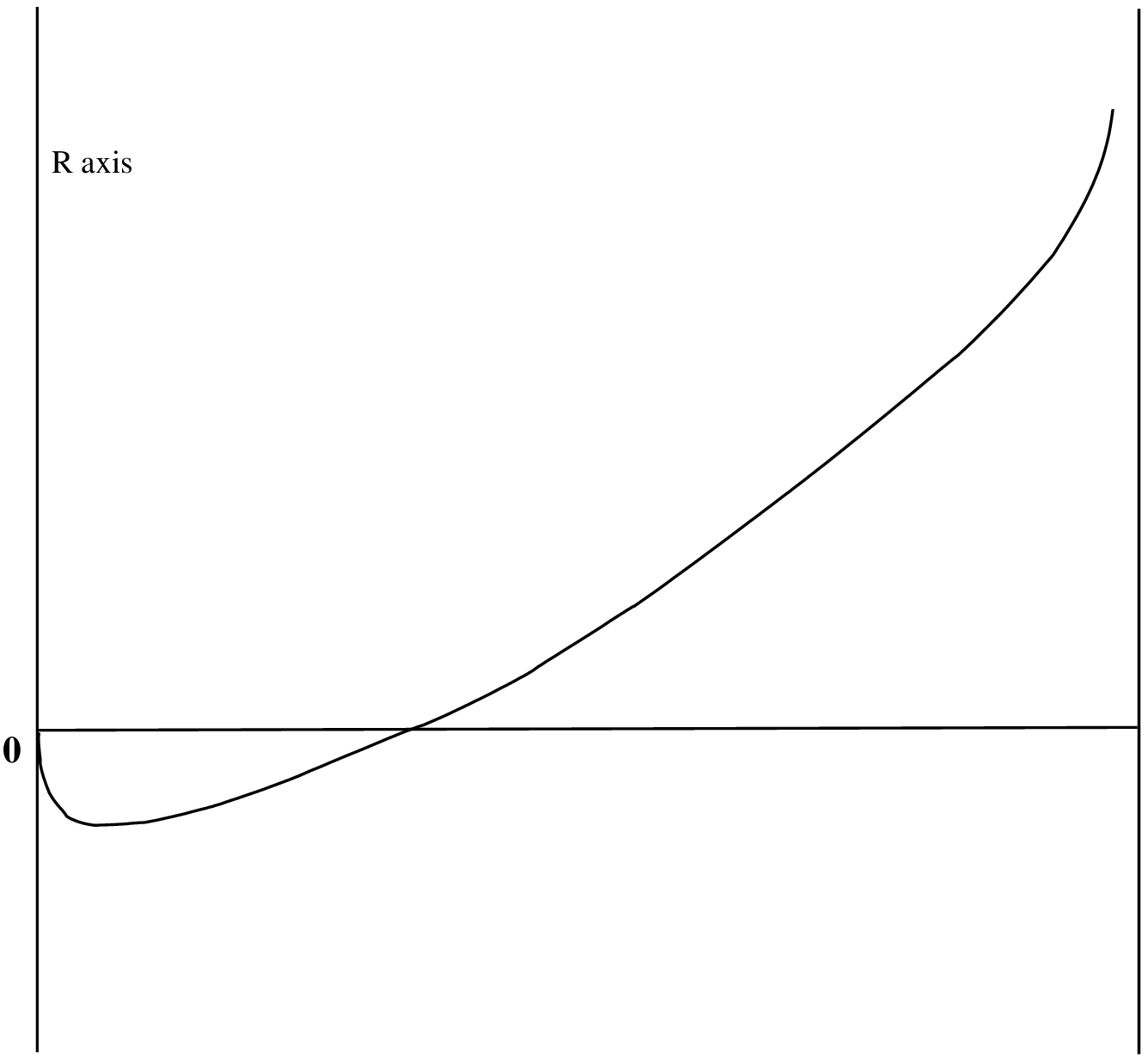}
    \caption[$n=2, H =\frac{1}{2}$, $d = - 1$]{$n=2, H =\frac{1}{2}$, \\ $d = - 1$}
    \label{F-tra-2d}
\end{minipage}\hfill
\end{center}
\end{figure}
\end{pb1-figs}\bigskip

\begin{pb1-figs}
\begin{figure}[h]
\begin{center}
\begin{minipage}[c]{4.5cm}
\includegraphics[width=4.5cm]{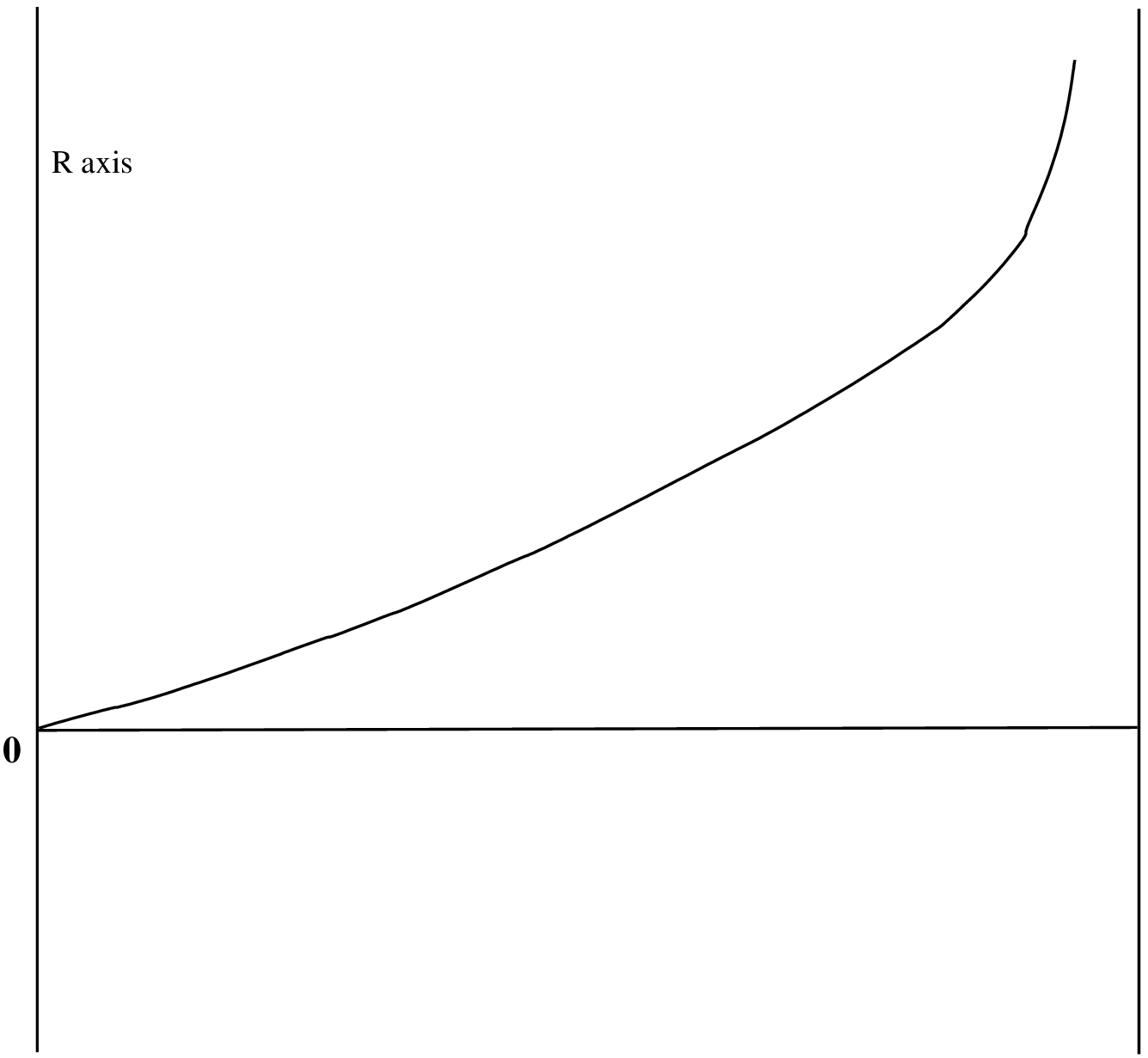}
    \caption[$n=2, H =\frac{1}{2}, 0 < d < 1$]{$n=2, H =\frac{1}{2}$, \\ $0 < d < 1$}
    \label{F-tra-2b}
\end{minipage}\hfill
\begin{minipage}[c]{4.5cm}
    \includegraphics[width=4.5cm]{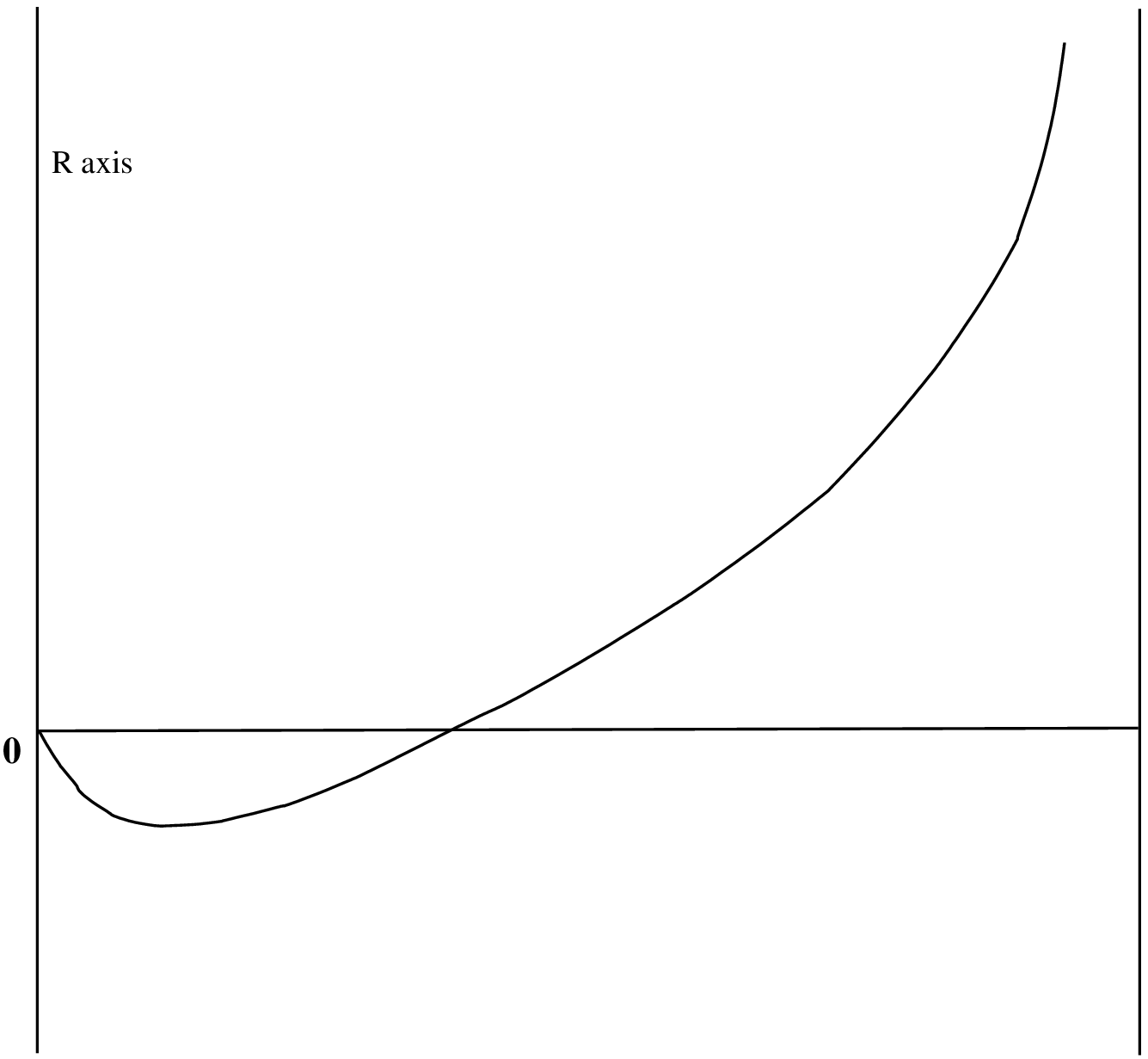}
    \caption[$n=2, H =\frac{1}{2}, -1 < d < 0$]{$n=2, H =\frac{1}{2}$, \\ $-1 < d < 0$}
    \label{F-tra-2e}
\end{minipage}\hfill
\end{center}
\end{figure}
\end{pb1-figs}


\newpage
\section{Applications, embedded minimal hypersurfaces with\\ boundary
contained in a slice}\label{S-appli}


In this section we give some results in which we use the
$H$-hypersurfaces constructed in Section \ref{S-examples} as
barriers. \bigskip

Recall from Section \ref{S-examples} that for $0 < H \le
\frac{n-1}{n}$ and for $d=0$, there exist simply-connected rotation
$H$-hypersurfaces $\cS_{H}$ which are entire vertical graphs going
to infinity at infinity. The unit normal of $\cS_{H}$ points
upwards. We call $\check{\cS}_{H}$ the symmetric of $\cS_{H}$ with
respect to the slice $\HH^n \times \{0\}$. Its unit normal points
downwards. We call $\cC(\cS_{H})$ the mean convex side of $\cS_{H}$
(\ie the connected component of the complement of $\cS_{H}$ into
which the unit normal points). We consider the set $\cR$ of
hypersurfaces obtained from $\cS_{H}$ and $\check{S}_{H}$ by
vertical or horizontal translations in $\HH^n \times \R$. We denote
by $\cC(\cS)$ the mean convex side of a hypersurface $\cS \in \cR$.
\bigskip

The following Proposition generalizes to higher dimensions the
\emph{convex hull lemma} given in \cite{NSST08}, Lemma 2.1.\bigskip

\begin{prop}[Convex hull lemma]\label{P-appl-1}
Given $K$ a compact subset in $\HH^n \times \R$, let $\cF_{K}^{H}$
denote the subset of domains $B$ in $\HH^n \times \R$ which contain
$K$ and such that $B = \cC(\cS)$ for some $\cS \in \cR$. Let $M$ be
a compact connected immersed hypersurface in $\HH^n \times \R$ with
mean curvature $H$.
\begin{enumerate}
    \item If $H$ is a constant in $]0, \frac{n-1}{n}]$, then $M
    \subset \cF_{\partial M}^{H}$.
    \item If $0 < H(x) \le \frac{n-1}{n}$ for all $x \in M$, then $M
    \subset \cF_{\partial M}^{(n-1)/n}$.
\end{enumerate}
\end{prop}\bigskip

\pf Because $M$ is compact, taking into account the asymptotic
behaviour of $\cS_{H}$ (see Theorems~\ref{T-h3-r0}-\ref{T-h3-rp}),
there exists some vertical translation $\tau$ such that $M \subset
\cC(\tau(\cS_{H}))$ so that the set of hypersurfaces in $\cR$ such
that $M \subset \cC(\cS)$ is non empty. Take any $\cS \in \cR$ such
that $M \subset \cC(\cS)$ and translate $\cS$ horizontally along
some geodesic until it touches $M$ at some point $p$. We claim that
$p$ cannot be an interior point. Indeed, assume that $p$ is an
interior point and let $p_0$ denote the projection of $p$ onto
$\HH^n$. Both hypersurfaces $\cS$ and $M$ would be vertical graphs
near $p_0$, corresponding respectively to functions $u, v$ such that
$u(p_0) = v(p_0)$ and $u \le v$ in a neighborhood of $p_0$. By the
maximum principle, this would imply that $M = \cS$ a contradiction.
The Proposition follows.

\hfill  \qed
\bigskip

In the applications below, we consider a hypersurface $\Gamma$ in
$\HH^n$ with the following properties.

\begin{equation}\label{E-appl-gam}
\left\{%
\begin{array}{ll}
    \Gamma & \text{is smooth, compact, connected, embedded,} \\
    \Gamma = \partial \Omega , & \Omega \text{ a bounded domain in } \HH^n,\\
    \Gamma & \text{has all its principal curvatures } > 1,\\
\end{array}%
\right.
\end{equation}

where the principal curvatures are taken with respect to the unit
normal to $\partial \Omega $ pointing inwards.\bigskip

Given a hypersurface $\Gamma$ satisfying Properties
(\ref{E-appl-gam}), there exists some radius $R$ such that for any
point $p$, the ball $B_{p,R} \subset \HH^n$ with radius $R$ is
tangent to $p$ at $\Gamma$ and $\Gamma \subset B_{p,R}$. We denote
by

\begin{equation}\label{E-appl-gam2}
\cS_{p,+} \text{ ~and~ } \cS_{p,-}
\end{equation}

the two hypersurfaces in $\cR$ passing through the sphere $\partial
B_{p,R}$ and symmetric with respect to the slice $\HH^n \times
\{0\}$.\bigskip

We first prove an existence result for a Dirichlet problem.
\bigskip

\begin{prop}\label{P-appl-2}
Let $\Omega \subset \HH^n \times \{0\}$ be a bounded domain with
smooth boundary $\Gamma$ satisfying (\ref{E-appl-gam}). Then, for
any $H, 0 < H \le \frac{n-1}{n}$, there exists a vertical graph
$M_{\Gamma}$ over $\Omega$ in $\HH^n \times \R$, with constant mean
curvature $H$ with respect to the upward pointing normal. This means
that there exists a function $u : \Omega \to \R$, smooth up to the
boundary, such that $u|_{\Gamma} = 0$, and whose graph
$\ens{(x,u(x))}{x \in \Omega}$ has constant mean curvature $H$ with
respect to the unit normal pointing upwards.
\end{prop}\bigskip

\textbf{Remark.} ~The graph $M_{\Gamma}$ having positive mean
curvature with respect to the upward pointing normal, must lie below
the slice $\HH^n \times \{0\}$. The symmetric $\check{M}_{\Gamma}$
with respect to the slice lies above the slice and has positive mean
curvature with respect to the normal pointing downwards. \bigskip

\textbf{Proof of Proposition \ref{P-appl-2}}\bigskip

\noib We first consider the case $H = \frac{n-1}{n}$.\bigskip

By our assumption on $\Gamma$, using the hypersurfaces
(\ref{E-appl-gam2}) and the Convex hull lemma, Proposition
\ref{P-appl-1}, any solution to our Dirichlet problem must be
contained in $\cC(\cS_{p,-}) \cap \cC(\cS_{p,+})$. This provides a
priori height estimates and boundary gradient estimates on the
solution. \bigskip

We could use \cite{Spr08} and classical elliptic theory \cite{GT83},
to get existence for our Dirichlet problem when $H=\frac{n-1}{n}$.
We shall instead apply \cite{Spr08} directly. Indeed, in our case,
the mean curvature $H_{\Gamma}$ of $\Gamma$ satisfies $H_{\Gamma} >
1 = H \frac{n}{n-1}$, and the Ricci curvature of $\HH^n$ satisfies
$\ric = - (n-1) \ge - \frac{n^2}{n-1}H^2$. Theorem 1.4 in
\cite{Spr08} states that under theses assumptions there exists a
vertical graph over $\Omega$ with boundary $\Gamma$ and constant
mean curvature $H=\frac{n-1}{n}$.\bigskip

\noib We now consider the case $0 < H \le \frac{n-1}{n}$. \bigskip

We use the graphs constructed previously as barriers to obtain a
priori height estimates and apply the interior and global gradient
estimates of \cite{Spr08} to conclude. \bigskip

We consider the Dirichlet problem $(P_t)$ for $0 \le t \le 1$,

\begin{equation*}\label{E-appl-5}
\left\{%
\begin{array}{cccc}
\mathrm{div}\big(\dfrac{\nabla u}{W}\big) & = &
t \, (n-1) & \text{in~ }\Omega\\[5pt]
u & = & 0 & \text{on~ } \Gamma \\
\end{array}%
\right.
\end{equation*}

where $u\in C^2(\Omega)$ is the height function, $\nabla u$ its
gradient and $W=(1 +|\nabla u|^2)^{1/2}$, and where the gradient and
the divergence are taken with respect to the metric on $\HH^n$. This
is the equation for \emph{vertical} $H$-graphs in $\HH^n \times \R$.
It is elliptic of divergence type. \bigskip

By the  first step, we have obtained the solution $u_1$ for the
Dirichlet problem $(P_1)$. The solution for $(P_0)$ is the trivial
solution $u_0=0$. By the maximum principle, using the fact that
vertical translations are positive isometries for the product
metric, and the existence of the solutions $u_1$ and $u_0$, we have
that any $C^1(\overline{\Omega})$ solution $u_t$ of the Dirichlet
problem $(P_t)$ stays above $u_1$ and below $u_0.$ This yields a
priori height and boundary  gradient estimates, independently of $t$
and $u_t$. Global gradient estimates follow Theorem 1.1 and Theorem
3.1 in \cite{Spr08}. We have therefore $C^1(\overline{\Omega})$ a
priori estimates independently of $t$ and $u_t$. The existence of
the solution $u_t$ for $0<t<1$ now follows from classical elliptic
theory, see \cite{GT83} or Theorem  A.7 in \cite{BaS98}. \bigskip

This completes the proof of Proposition \ref{P-appl-2}. \hfill  \qed
\bigskip

We now generalize to higher dimensions results obtained in
\cite{NSST08}.\bigskip

\newpage

\begin{thm}\label{T-appl-6}
Let $M$ be an embedded compact connected $H$-hypersurface in $\HH^n
\times \R$, with \goodbreak $0 < H \le \frac{n-1}{n}$. Assume that
the boundary $\Gamma$ is an $(n-1)$-submanifold in $\HH^n \times
\{0\}$ satisfying (\ref{E-appl-gam}).
\begin{enumerate}
    \item The hypersurface $M$ is either the graph $M_{\Gamma}$
    given by Proposition \ref{P-appl-2} or its symmetric
    $\check{M}_{\Gamma}$.
    \item Assume furthermore that $\Gamma$ is symmetric with respect
    to some hyperbolic hyperplane $P$ in $\HH^n \times \{0\}$ and
    that each connected component of $\Gamma \setminus P$ is a graph
    above $P$. Then $M$ is symmetric with respect to the vertical
    hyperplane $P \times \R$ and each connected component of $M
    \setminus P\times \R$ is a horizontal graph. In particular, if
    $\Gamma$ is an $(n-1)$-sphere, the hypersurface $M$ is part of
    the rotation surface given by Theorem \ref{T-h3-r0}.
\end{enumerate}
\end{thm}\bigskip

\textbf{Proof of Theorem \ref{T-appl-6}.}\bigskip

\noi Let $\Omega$ be the bounded domain such that $\Gamma =
\partial \Omega$ and let $\cC = \overline{\Omega} \times \R$ be the vertical
cylinder above $\overline{\Omega}$. We claim that $M \subset \cC$
and that $M \cap \cC = \Gamma$. Indeed, at each $p \in \Gamma$, we
have the hypersurfaces $\cS_{p,+}$ and $\cS_{p,-}$ given by
(\ref{E-appl-gam2}). It follows from the Convex hull lemma,
Proposition \ref{P-appl-1}, that $M$ is in the convex hull of such
hypersurfaces and hence that $M \subset \cC$ and $M \cap \cC =
\Gamma$.\bigskip

\noi By Proposition \ref{P-appl-2}, we have two vertical graphs
above $\Omega$, $M_{+} \subset \HH^n \times \R_{+}$ with constant
mean curvature $H$ with respect to the normal pointing downwards and
$M_{-} \subset \HH^n \times \R_{-}$ with constant mean curvature $H$
with respect to the normal pointing upwards. \bigskip

\noi We claim that $M$ is a vertical graph contained either in
$\HH^n \times \R_{+}$ or in $\HH^n \times \R_{-}$. If not, making
reflexions with respect to slices $\HH^n \times \{t\}$ starting from
$t_{+}$ the highest height on $M$ we would obtain a contradiction by
the maximum principle. If $M$ were not contained in one of the
half-spaces, we would have highest and lowest interior points at
which the normal would point downwards, \resp upwards by the maximum
principle.


\bigskip

\noi We claim that $M=M_{+}$ or $M=M_{-}$. Assume that $M \subset
\HH^n \times \R_+$ (the proof is similar if $M$ is contained in the
lower half-space). Translating $M_{+}$ vertically upwards very far
and then coming down, we see that $\tau(M_{+})$ cannot touch $M$
before the boundaries coincide (maximum principle). It follows that
$M$ must be below $M_{+}$. Doing the same thing with $M$, we see
that $M$ must be above $M_{+}$. It follows finally that $M=M_{+}$.

\noi Assume now that $\Gamma$ is symmetric with respect to a
hyperbolic hyperplane $P$ and assume that each connected component
of $\Gamma\setminus P$ is a horizontal graph. We can then use
Alexandrov Reflection Principle in vertical hyperplanes $P_t\times
\R$ in ambient space, obtained by applying horizontal translations
along geodesics orthogonal to $P$, to the vertical hyperplane
$P\times \R$ of symmetry of $\Gamma$, and conclude that $M$ is
symmetric with respect to $P\times \R$. Moreover, Alexandrov
Reflection Principle ensures that each connected component of
$M\setminus P\times \R$ is a horizontal graph.\bigskip

When $\Gamma$ is an $(n-1)$-sphere,  we can apply the preceding
result to prove that $M$ is rotationally symmetric.

\hfill \qed
\bigskip

Recall from \cite{BS08a} that the height of the family of minimal
catenoids in $\HH^n \times \R$ is $\frac{\pi}{n-1}$.\bigskip

\begin{thm}\label{T-appl-7}
Let $\Gamma$ satisfy (\ref{E-appl-gam}). Consider two copies of
$\Gamma$ in different slices $\Gamma_{+} = \Gamma \times \{a\}$ and
$\Gamma_{-} = \Gamma \times \{-a\}$ for some $a>0$. Let $M$ be a
compact connected embedded $H$-hypersurface such that $\partial M =
\Gamma_{+} \cup \Gamma_{-}$, with $0 < H \le \frac{n-1}{n}$. Assume
that $2a \ge \frac{\pi}{n-1}$.
\begin{enumerate}
    \item Assume that $\Gamma$ is symmetric with respect to a hyperbolic
    hyperplane $P$ and that each connected component of $\Gamma\setminus P$
    is a graph above $P$. Then $M$ is symmetric with respect to the vertical
    hyperplane  $P\times \R$ and each connected component of $M\setminus
    P\times \R$ is a horizontal graph.
    \item Assume that $\Gamma$ is an $(n-1)$-sphere. Then $M$ is part of
    the complete embedded rotation hypersurface given by
    Theorem~\ref{T-h3-r0} and \ref{T-h3-rp} and containing $\Gamma$.
    It follows that $M$ is symmetric with respect to the slice $\HH^n \times
    \{0\}$ and the parts of $M$ above and below the slice of symmetry are
    vertical graphs.
\end{enumerate}
\end{thm}\bigskip

\textbf{Proof of Theorem \ref{T-appl-7}}.\bigskip

\noi Let $\Omega_{+} = \Omega \times \{a\}$ and $\Omega_{-} = \Omega
\times \{-a\}$. By the Convex hull Lemma, Proposition
\ref{P-appl-1}, using the hypersurfaces given by (\ref{E-appl-gam2})
we have that $M \cap \overline{\mathrm{ext}(\Omega_{+})} = \Gamma_+$
and $M \cap \overline{\mathrm{ext}(\Omega_{-})} = \Gamma_-$.\bigskip

\noi We claim that $M \cap (\overline{\Omega} \times \R) = \Gamma_+
\cap \Gamma_-$. Let $M_{\Gamma,a}$ be the graph above
$\overline{\Omega_+}$ contained in $\HH^n \times [a, \infty[$ and
$M_{\Gamma,-a}$ be the graph below $\overline{\Omega_-}$ contained
in $\HH^n \times ]-\infty, a]$, given by Theorem~\ref{T-appl-6}.
\bigskip

Consider $\widetilde{M} = M_{\Gamma,a} \cap M \cap M_{\Gamma,-a}$
oriented by the mean curvature vector of $M$ by continuity. Take the
family of (minimal) catenoids symmetric with respect to $\HH^n
\times \{0\}$ with rotation axis some $\{\bullet \} \times \R$.
Coming from infinity with such catenoids, using the assumption that
$2a \ge \frac{\pi}{n-1}$ and the fact that the catenoids have height
$< \frac{\pi}{n-1}$, we see that one catenoid will eventually touch
$\widetilde{M}$ at some interior point in $M$. This implies that the
normal to $M$ at this point is the same as the normal to the
catenoid at the same point (maximum principle) and hence that the
normal to $M$ points inside $\widetilde{M}$.\bigskip

Assume that $M \cap (\Omega \times \{a\}) \not = \emptyset$ (\resp
that $M \cap (\Omega \times \{-a\}) \not = \emptyset$). Then at the
highest point of $M$ the normal would be pointing upwards (\resp
downwards) and we would get a contradiction with the maximum
principle by considering the horizontal slice (a minimal
hypersurface) at this point.\bigskip

Finally, $M \cap (\overline{\Omega} \times \R) = \Gamma_+ \cap
\Gamma_-$ and the normal to $M$ points inside $M \cup \Omega_+ \cup
\Omega_-$.\bigskip

\noi To conclude, we use Alexandrov Reflection Principle in vertical
hyperplanes $P_t\times \R$ in ambient space, obtained by applying
horizontal translations along the horizontal geodesic orthogonal to
$P,$  to the hyperplane $P\times \mathbb R$ of symmetry of $\Gamma$.
We conclude that $M$ is symmetric about $P\times \R$ and  that each
connected component of $M\setminus P\times \R$  is a horizontal
graph. This complete the proof of the first statement in the
theorem. \bigskip

If $\Gamma$ is spherical then $M$ is a rotation hypersurface. As the
mean curvature vector points into the region of ambient space that
contains the axes, by the geometric classification of the rotation
$H$-hypersurfaces with constant mean curvature $H \le (n-1)/n$ given
by Theorems~\ref{T-h3-r0} and \ref{T-h3-rp}, it follows that $M$ is
part of a complete embedded rotation hypersurface $\overline{M}$. It
follows that $\overline{M}$ has a slice of symmetry at $\HH^n \times
\{0\}$ and each connected component of  $\overline{M}$  above and
below $t=0$ is a complete vertical graph  over the exterior of a
round ball in $t=0$.

\hfill \qed \bigskip

\appendix

\bigskip
\section{Vertical flux formula in $\HH^n \times \R$}\label{S-vflux}


Let $f : M \looparrowright \Mh := \HH^n \times \R$ be an isometric
immersion. Let $h$ denote the function $h : \Mh \rightarrow \R$,
such that $h(x,t) = t$ and let $h_{M} = h|_{M}$ be the restriction
of the function $h$ to the hypersurface $M$, \ie the height function
of $M$. We let $\gh$ be the (product) metric on $\Mh$ and $\Delta_M$
be the (non-positive) Laplacian on $M$, for the induced metric $g :=
f^{*}\gh$. \bigskip

\begin{prop}\label{P-a-flux-1}
With the above notations we have
$$\Delta_M h_M = n \gh(\overrightarrow{H},\partial_t)$$
where $\overrightarrow{H}$ is the (normalized) mean curvature vector
of the immersion and $\partial_t$ the vertical vector-field along
$\R$.
\end{prop}\bigskip

\textbf{Remark.}~ When $f$ admits a unit normal field $N_M$ the
above formula boils down to $\Delta_M h_M = n H v_M$ where $H$ is
the (normalized) mean curvature in the direction $N_M$ and $v_M$ the
vertical component of $N_M$, $v_M := \gh(N_M,\partial_t)$. \bigskip

\pf Take a local orthonormal frame $\{E_i\}_{i=1}^n$ for $M$ near a
point $m \in M$ and extend it locally in a neighborhood of $m$ in
$\Mh$. Then

\begin{eqnarray*}
  \Delta_M h_M &=& \sum_{i=1}^n \big\{ \big(E_i \cdot (E_i\cdot h_M) \big)
  -(D_{E_i}E_i) \cdot h_M \big\} \\
   &=& \sum_{i=1}^n \big\{ E_i \cdot \big( dh_M(E_i) \big)
  -dh_M(D_{E_i}E_i) \big\} \\
   &=& \sum_{i=1}^n \big\{ E_i \cdot \big( dh(E_i) \big)
  -dh(D_{E_i}E_i) \big\}\\
  &=& \sum_{i=1}^n \big\{(\Dh_{E_i}dh)(E_i) + dh(\Dh_{E_i}E_i) - dh(D_{E_i}E_i)
  \big\}.
\end{eqnarray*}\bigskip

In the product space $\Mh = \HH^n \times \R$, we have $\Dh_{E}dh=0$
for all $E \in \cX(\Mh)$. It follows that

$$\Delta_M h_M = \sum_{i=1}^n dh(\Dh_{E_i}E_i - D_{E_i}E_i) =
\sum_{i=1}^n dh(A(E_i,E_i))$$

where $A$ is the second fundamental form of the immersion. Finally,

$$\Delta_M h_M = dh(\mathrm{Tr}(A)) = n dh(\overrightarrow{H})$$

which is the formula in the Theorem.

\hfill \qed\bigskip

\begin{cor}\label{C-a-flux-2}
Let $\Omega$ be a compact domain on $M$ with unit inner normal
$\nu_{\partial \Omega}$ to $\partial \Omega$ in $\Omega$. Then
$$\int_{\Omega} \Delta_M h_M \, d\mu_M = - \int_{\partial \Omega}
dh_M(\nu_{\partial \Omega}) \, d\sigma_{\partial \Omega} = n
\int_{\Omega} \gh(\overrightarrow{H},\partial_t) \, d\mu_M.$$
\end{cor}\bigskip

\pf Divergence Theorem.

\hfill \qed \bigskip

\textbf{Applications to rotation $H$-hypersurfaces}\bigskip

Let us consider a rotation hypersurface $M$ given by the
parametrization

$$X(\rho,\xi) = \big( \tanh(\rho /2)\xi , \lambda(\rho) \big)$$

with $\rho > 0$ and $\xi \in S^{n-1}$ and choose the unit normal
pointing upwards.
\bigskip

Consider the domain

$$\Omega (\rho_0,\rho) := X([\rho_0,\rho]\times S^{n-1}) \subset M.$$

We have

$$X_{\rho}(\rho,\xi) \big( \frac{\xi}{2 \cosh^2(\rho/2)}, \dot{\lambda}(\rho) \big),$$

$$v_M(\rho,\xi) = (1+\dot{\lambda}^2(\rho))^{-1/2},$$

$$d\mu_M = (1+\dot{\lambda}^2(\rho))^{1/2} \sinh^{n-1}(\rho) \, d\rho \, d\mu_S,$$

$$\nu_{\partial \Omega (\rho_0,\rho)(X(\rho,\xi))} = (1+\dot{\lambda}^2(\rho))^{-1/2}
X_{\rho}(\rho,\xi),$$

$$d\sigma_{X(\{\rho\} \times S^{n-1})} = \sinh^{n-1}(\rho) \, d\mu_{S}.$$

\bigskip

The above Corollary applied to $\Omega (\rho_0,\rho)$ gives

$$- \mathrm{Vol}(S^{n-1}) \sinh^{n-1}(t) \dot{\lambda}(t) (1+\dot{\lambda}^2(t))^{-1/2}
\Big|_{\rho_0}^{\rho} = - n \mathrm{Vol}(S^{n-1})
\int_{\rho_0}^{\rho} H(t) \sinh^{n-1}(t) \, dt.$$

Looking for rotation surfaces with \emph{constant mean curvature}
$H$ we find

$$\sinh^{n-1}(\rho) \dot{\lambda}(\rho) (1+\dot{\lambda}^2(\rho))^{-1/2}
= n H \int_{\rho_0}^{\rho} \sinh^{n-1}(t) \, dt + F(\rho_0)$$

where the constant $F(\rho_0) := \sinh^{n-1}(\rho_0)
\dot{\lambda}(\rho_0) (1+\dot{\lambda}^2(\rho_0))^{-1/2}$ is the
flux through $X(\{\rho_0\} \times S^{n-1})$.

\vspace{1.5cm}

\def\refname{References}

\vspace{1.5cm}

\begin{flushleft}
\begin{tabular}{lll}
Pierre B\'{e}rard &&  Ricardo Sa Earp\\
Universit\'{e} Joseph Fourier && Departamento de Matem\'{a}tica\\
Institut Fourier - Math\'{e}matiques (UJF-CNRS) && PUC-Rio \\
 B.P. 74 && 22453-900, Rio de Janeiro - RJ - Brazil\\
38402 Saint Martin d'H\`{e}res Cedex - France\\
{\em E-mail}: Pierre.Berard@ujf-grenoble.fr &&  {\em E-mail}: earp@mat.puc-rio.br
\end{tabular}
\end{flushleft}

\end{document}

%% file: macros.tex
\newcommand{\noi}{\noindent }
\newcommand{\noib}{\noindent $\bullet $~}
\newcommand{\pointir}{\discretionary{.}{}{.\kern.35em---\kern.7em}}
\newcommand{\ie}{\emph{i.e. }}

\newcommand{\resp}{\emph{resp. }}

\usepackage{latexsym}
\usepackage{amsmath}
\usepackage{amssymb}

\newcommand{\Mh}{\widehat{M}}

\newcommand{\gh}{\widehat{g}}
\newcommand{\Dh}{\widehat{D}}

\newcommand{\ric}{\mathrm{Ric}}

\def\buildo#1\over#2{\mathrel{\mathop{\null#2}\limits^{#1}}}
\def\buildu#1\under#2{\mathrel{\mathop{\null#2}\limits_{#1}}}

\newcommand{\B}{\mathbb{B}}

\newcommand{\R}{\mathbb{R}}

\newcommand{\HH}{\mathbb{H}}
\newcommand{\PP}{\mathbb{P}}
\newcommand{\V}{\mathbb{V}}

\newcommand{\cC}{{\mathcal C}}
\newcommand{\cD}{{\mathcal D}}

\newcommand{\cF}{{\mathcal F}}

\newcommand{\cK}{{\mathcal K}}

\newcommand{\cN}{{\mathcal N}}

\newcommand{\cR}{{\mathcal R}}
\newcommand{\cS}{{\mathcal S}}
\newcommand{\cT}{{\mathcal T}}
\newcommand{\cU}{{\mathcal U}}

\newcommand{\cX}{{\mathcal X}}

\newcommand{\ens}[2]{\{ #1 ~|~ #2 \}}

\newcommand{\pf}{\par{\noindent\textbf{Proof.~}}}

\newbox\qedbox
\setbox\qedbox=\hbox{$\Box$}
\newcommand{\qed}{\hfill\penalty10000\copy\qedbox\par\medskip}
\def\refname{References}

\newtheorem{thm}{Theorem}[section]
\newtheorem{thm-def}[thm]{Theorem and Definition}

\newtheorem{prop}[thm]{Proposition}

\newtheorem{prop-def}[thm]{Proposition and Definition}

\newtheorem{cor}[thm]{Corollary}